\documentclass{mcom-l}

\usepackage{amssymb}

\usepackage{graphicx}

\newtheorem{theorem}{Theorem}[section]
\newtheorem{lemma}[theorem]{Lemma}

\theoremstyle{definition}
\newtheorem{definition}[theorem]{Definition}

\theoremstyle{remark}
\newtheorem{remark}[theorem]{Remark}

\numberwithin{equation}{section}

\begin{document}

\title{Localization of Elliptic Multiscale Problems}


\author{Axel M\r{a}lqvist}
\address{Department of Information Technology, Uppsala University, Box 337, SE-751 05 Uppsala, Sweden}
\curraddr{}
\email{}
\thanks{A. M\r{a}lqvist is supported by The G\"{o}ran Gustafsson Foundation and The Swedish Research Council.}

\author{Daniel Peterseim}
\address{Rheinische Friedrich-Wilhelms-Universit\"at Bonn, Institute for Numerical Simulation, Wegelerstr. 6, 53115 Bonn, Germany}
\curraddr{}
\email{peterseim@ins.uni-bonn.de}
\thanks{D. Peterseim was supported by the Humboldt-Universt\"at zu Berlin and the DFG Research Center Matheon Berlin through project C33. }


\date{}

\dedicatory{}

\begin{abstract}
This paper constructs a local generalized finite element basis for elliptic problems with heterogeneous and highly varying coefficients. The basis functions are solutions of local problems on vertex patches. The error of the corresponding generalized finite element method decays exponentially with respect to the number of layers of elements in the patches. Hence, on a uniform mesh of size $H$, patches of diameter $H\log(1/H)$ are sufficient to preserve a linear rate of convergence in $H$ without pre-asymptotic or resonance  effects.
The analysis does not rely on regularity of the solution or scale separation in the coefficients.
This result motivates new and justifies old classes of variational multiscale methods.
\end{abstract}

\maketitle

\newcommand{\tri}{\mathcal{T}}
\newcommand{\triH}{\tri_H}
\newcommand{\trih}{\tri_h}
\newcommand{\grad}{\nabla}
\newcommand{\dist}{\operatorname*{dist}}
\newcommand{\diam}{\operatorname*{diam}}
\newcommand{\support}{\operatorname*{supp}}
\newcommand{\dimension}{\operatorname*{dim}}
\newcommand{\conv}{\operatorname*{conv}}
\newcommand{\inter}{\operatorname*{int}}
\newcommand{\Inodal}{\mathfrak{I}_{H}}
\newcommand{\Inodalh}{\mathfrak{I}_{h}}
\newcommand{\Tfine}{\mathfrak{F}}
\newcommand{\ums}{u^{\operatorname*{ms}}_H}
\newcommand{\vms}{u^{\operatorname*{ms}}_H}
\newcommand{\umsloc}[1]{u^{\operatorname*{ms}}_{H,#1}}
\newcommand{\umsloch}[1]{u^{\operatorname*{ms},h}_{H,#1}}
\newcommand{\zmsloch}[1]{z^{\operatorname*{ms},h}_{H,#1}}
\newcommand{\umsloctilde}[1]{\tilde{u}^{\operatorname*{ms}}_{H,#1}}
\newcommand{\ufine}{u^{\operatorname*{f}}}
\newcommand{\ufineh}{u^{\operatorname*{f}}_h}
\newcommand{\VmsHh}{V^{\operatorname*{ms},h}_{H}}
\newcommand{\umsh}{u^{\operatorname*{ms},h}_H}
\newcommand{\Creg}{C_{A,\Omega}}
\newcommand{\Cint}{C_{\Inodal}}
\newcommand{\Cinv}{C_{\Inodal}'}
\newcommand{\Cinverse}{C_{\operatorname*{inv}}}
\newcommand{\Col}{C_{\operatorname*{ol}}}
\newcommand{\Cov}{C_{\operatorname*{ov}}}
\newcommand{\Cco}{C_{\operatorname*{co}}}
\newcommand{\VmsH}{V^{\operatorname*{ms}}_H}
\newcommand{\VmsHk}{V^{\operatorname*{ms}}_{H,k}}
\newcommand{\VmsHkh}{V^{\operatorname*{ms},h}_{H,k}}
\newcommand{\Vf}{V^{\operatorname*{f}}}
\newcommand{\vf}{v^{\operatorname*{f}}}
\newcommand{\Vhxf}{V_{h}^{\operatorname*{f}}}
\newcommand{\tnorm}[1]{\left|\left|\left|#1\right|\right|\right|}
\newcommand{\Tnorm}[1]{\biggl|\biggl|\biggl|#1\biggr|\biggr|\biggr|}

\section{Introduction}
This paper considers the numerical solution of second order elliptic problems with strongly heterogeneous and highly varying (non-periodic) coefficients. The heterogeneities and oscillations of the coefficients may appear on several non-separated scales.
It is well known that classical polynomial based finite element methods perform arbitrarily badly for such problems, see e.g.~\cite{MR1648351}. To overcome this lack of performance, many methods that are based on general (non-polynomial) ansatz functions have been developed. Early works \cite{MR701094,MR1286212}, that essentially apply to one dimensional problems, have been generalized to the multi-dimensional case in several ways during the last fifteen years, see e.g.~\cite{MR1660141,MR1455261,resbub}. In these methods the problem is split into coarse and (possibly several) fine scales. The fine scale effect on the coarse scale is either computed numerically or modeled analytically. The resulting modified coarse problem can then be solved numerically and its solution contains crucial information from the fine scales. Although many of these approaches show promising results in practice, their convergence analysis usually assumes certain periodicity and scale separation.

For problems with general $L^{\infty}$ coefficients, the paper \cite{BabLip11} gives error bounds for a generalized finite element method that involves the solutions of local eigenvalue problems. The construction in \cite{MR2721592,OwhadiZ11} depends only on the solution of the original problem on certain subdomains. However, the size of these subdomains strongly depends on the mesh size. This dependence is suboptimal with respect to the theoretical statement given in \cite{ALBasis1}, that is, for any shape regular mesh of size $H$ there exist $\mathcal{O}\left(  \left(\log(1/H)\right)^{d+1}\right)$ local (non-polynomial) basis functions per nodal point such that the error of the corresponding Galerkin solution $u_{H}$ satisfies the estimate $\left\Vert u-u_{H}\right\Vert _{H^{1}(\Omega) }\leq C_{g}H$ with a constant $C_{g}$ that depends on the right-hand side $g$ and the global bounds of the
diffusion coefficient but not on its variations. The derivation in \cite{ALBasis1} is not constructive in the sense that it involves the solution of the (global) original problem with specific right hand sides.

In this paper, we show that such a (quasi-)optimal basis can indeed be constructed by solving only local problems on element patches. We use a modified nodal basis similar to the one presented in \cite{MR2319044} and prove that these basis functions decay exponentially away from the node they are associated with. This exponential decay justifies an approximation using localized patches.

The precise setting of the paper is as follows. Let $\Omega\subset\mathbb{R}^{d}$ be a bounded Lipschitz domain with polygonal boundary
and let the diffusion matrix $A\in L^\infty\left(\Omega,\mathbb{R}%
_{\mathrm{sym}}^{d\times d}\right)  $ be uniformly elliptic:
\begin{equation}%
\begin{array}
[c]{r}%
0<\alpha( A,\Omega) :=\underset{x\in\Omega}{\operatorname{ess}\inf}%
\inf\limits_{v\in\mathbb{R}^{d}\setminus\{0\}}\dfrac{(A( x) v)\cdot v
}{v\cdot v},\\
\\
\infty>\beta( A,\Omega) :=\underset{x\in\Omega}{\operatorname{ess}\sup}%
\sup\limits_{v\in\mathbb{R}^{d}\setminus\{0\}}\dfrac{(A( x) v)\cdot v
}{v\cdot v}.
\end{array}
\label{defalphabeta}%
\end{equation}
Given $g\in L^{2}( \Omega) $, we seek $u\in V:=H_{0}^{1}( \Omega) $
such that%
\begin{equation}\label{e:model}
a\left(  u,v\right):=\int_{\Omega}(A\nabla u)\cdot \nabla
v =\int_{\Omega}gv=:G( v) \quad\text{for all } v\in V.
\end{equation}
The bilinear form $a$ is symmetric, coercive, bounded, and hence, \eqref{e:model} has a unique solution.

The main result of this paper (cf. Theorem \ref{thm3}) shows that the error $u-\umsloc{k}$ of the generalized finite element method, which is based on our new (local) basis functions mentioned above, is bounded as follows
\begin{equation*}
 \|A^{1/2}\nabla(u-\umsloc{k})\|_{L^2(\Omega)}\leq C_g H;
\end{equation*}
$H$ being the mesh size of the underlying coarse finite element mesh and $k\approx\log(1/H)$ referring to the number of layers of coarse elements that form the support of the localized basis functions. This estimate shows that our new numerical upscaling procedure is reliable  beyond strong assumptions like periodicity and scale separation.  Moreover, the error bound is stable with respect to perturbations arising from the discretization of the local problems. These results give a theoretical foundation for numerous previous experiments where exponential decay of a similar modified basis have been noticed, see e.g.~\cite{MRX}.

The outline of the paper is as follows. In Section~\ref{s:locbasis}, we derive a set of local basis functions and define the corresponding multiscale finite element method. The error analysis is done in Section~\ref{s:error}. Section~\ref{s:disc} is devoted to the discretization of the local problems. Section~\ref{s:numexp} presents numerical experiments, and Section~\ref{s:application} discusses the application of this theory to state-of-the-art multiscale methods.

\section{Local Basis}\label{s:locbasis}
In this section, we design a set of local basis functions for the multiscale problem under consideration. The construction is based on a regular (in the sense of \cite{CiarletPb}) finite element mesh $\triH$ of $\Omega$ into closed triangles ($d=2$) or tetrahedra ($d=3$). Subsection \ref{ss:classical} recalls the classical nodal basis with respect to $\triH$ and demonstrates its lack of approximation properties. Subsection \ref{ss:intpol} introduces a quasi interpolation operator used in the construction of the new basis. Subsection \ref{ss:multiscale} defines a modified (coefficient dependent) nodal basis and analyzes its approximation properties. This basis is then localized in Subsection \ref{ss:localization}.

\subsection{Classical Nodal Basis}\label{ss:classical}
Let $H:\overline\Omega\rightarrow\mathbb{R}_{>0}$ denote the $\triH$-piecewise constant mesh size function with $H\vert_T=\diam(T)=:H_T$ for all $T\in\triH$. The mesh size may vary in space. In practical applications, the mesh $\triH$ (resp. its size $H$) shall be determined by the accuracy which is desired or the computational capacity that is available but \emph{not} by the scales of the coefficients.

The classical (conforming) $P_1$ finite element space is given by
\begin{equation}\label{e:couranta}
S_H:=\{v\in C^0(\bar{\Omega})\;\vert\;\forall T\in\triH,v\vert_T \text{ is a polynomial of total degree}\leq 1\}.
\end{equation}
Let $V_H:=S_H\cap V$ denote the space of finite element functions that match the homogeneous Dirichlet boundary conditions. Let $\mathcal{N}$ denote the set of interior vertices of $\triH$. For every vertex $x\in\mathcal{N}$, let $\lambda_x\in S_H$ denote the corresponding nodal basis function (tent function), i.e.,
\begin{equation*}
 \lambda_x(x)=1\;\text{ and }\;\lambda_x(y)=0\quad\text{for all }y\neq x\in\mathcal{N}.
\end{equation*}
These nodal basis functions form a basis of $V_H$. The availability of such a local basis is a key property of any finite element method and ensures that the resulting system of linear algebraic equations is sparse.

The (unique) Galerkin approximation $u_H\in V_H$ satisfies
\begin{equation}\label{e:courantb}
a(u_H,v) = G(v) \quad\text{ for all }v\in V_H.
\end{equation}
The above method \eqref{e:courantb} is optimal with respect to the energy norm $\tnorm{\cdot}:=\tnorm{\cdot}_{\Omega}:=\|A^{1/2}\grad\cdot\|_{L^2(\Omega)}$ on $V$ which is induced by $a$,
\begin{equation}\label{e:opt}
 \tnorm{u-u_H}=\min_{v_H\in V_H}\tnorm{u-v_H}.
\end{equation}
Assuming that the solution $u$ is smooth, the combination of \eqref{e:opt} and standard interpolation error estimates yields the standard a priori error estimate
\begin{equation*}\label{e:apriori}
 \tnorm{u-u_H}\leq C\|H\|_{L^\infty(\Omega)}\|\nabla^2 u\|_{L^2(\Omega)}.
\end{equation*}
This estimate states linear convergence of the classical finite element method \eqref{e:courantb} as the maximal mesh width tends to zero. However, the regularity assumption is not realistic for the problem class under consideration. Moreover, even if the coefficient is smooth, it may oscillate rapidly, say at frequency $\varepsilon^{-1}$ for some small parameter $\varepsilon$. In this case, the asymptotic result is useless because $\nabla^2 u$ may oscillate at the same scale, a fact that is hidden in the constant $\|\nabla^2u\|_{L^2(\Omega)}\approx \varepsilon^{-1}$.
Unless $H\lesssim \varepsilon$, the above finite element space is unable to capture the behavior of the solution neither on the microscopic nor on the macroscopic level.
In what follows, we present a new method that resolves this issue.

\subsection{Quasi Interpolation}\label{ss:intpol}
The key tool in our construction will be some bounded linear surjective (quasi-) interpolation operator $\Inodal: V\rightarrow V_H$. The choice of this operator is not unique and a different choice might lead to a different multiscale method.
We have in mind the following modification of Cl\'ement's interpolation \cite{MR0400739} which is presented and analyzed in \cite[Section 6]{MR1706735}. Given $v\in V$, $\Inodal v := \sum_{x\in\mathcal{N}}(\Inodal v)(x)\lambda_x$ defines a (weighted) Cl\'ement interpolant with nodal values
\begin{equation}\label{e:clement}\textstyle
 (\Inodal v)(x):=\left(\int_{\Omega} v \lambda_x\right)\left/\left(\int_{\Omega} \lambda_x\right)\right.
\end{equation}
for $x\in\mathcal{N}$. The nodal values are weighted averages of the function over nodal patches $\omega_x:=\support\lambda_x$. Since the summation is taken only with respect to interior vertices $N$, this operator matches homogeneous Dirichlet boundary conditions.
\begin{subequations}\label{e:quasiint}
\end{subequations}

Recall the (local) approximation and stability properties of the interpolation operator $\Inodal$ \cite[Section 6]{MR1706735}: There exists a generic  constant $\Cint$ such that for all $v\in V$ and for all $T\in\triH$ it holds
\begin{equation}\label{e:interr}
 H_T^{-1}\|v-\Inodal v\|_{L^{2}(T)}+\|\nabla(v-\Inodal v)\|_{L^{2}(T)}\leq \Cint \| \nabla v\|_{L^2(\omega_T)},
\tag{\ref{e:quasiint}.a}
\end{equation}
where $\omega_T:=\cup\{K\in\triH\;\vert\;T\cap K\neq\emptyset\}$. The constant $\Cint$ depends on the shape regularity parameter $\rho$ of the finite element mesh $\triH$ (see \eqref{e:shapereg} below) but not on $H_T$.

Note that the above interpolation operator is not a projection, i.e., $v_H\in V_H$ does not equal its interpolation $\Inodal v_H$ in general.
However, the particular choice gives rise to the following lemma.
\begin{lemma}\label{l:invert} There exists a generic constant $\Cinv$ which only depends on $\rho$ but not on the local mesh size $H$, such that for all $v_H\in V_H$ there exists $v\in V$ with the properties
\begin{equation}\label{e:invert}
 \Inodal(v)=v_H,\quad \mbox{   }\quad\|\nabla v\|\leq \Cinv\|\nabla v_H\|,\quad \text{and}\quad \support v\subset\support v_H.
\tag{\ref{e:quasiint}.b}
\end{equation}
\end{lemma}
\begin{proof}
For every nodal basis function $\lambda_x$, $x\in\mathcal{N}$, there is some $b_x\in H^1_0(\omega_x)$ such that $\Inodal(b_x)=\lambda_x$ and $\|\nabla b_x\|\leq \Cinv'\|\nabla \lambda_x\|$ with some constant $\Cinv'$ that does not depend on $x$ and $H$. E.g., $b_x$ may be chosen as a standard cubic element bubble on an arbitrary element $T\subset\omega_x$ or a quadratic edge/face bubble related to an arbitrary edge/face of $\triH$ interior to $\omega_x$. One might as well choose $b_x$ to be nodal interpolation of those bubbles in a finite element space that correponds to some uniform refinement of $\triH$.

Given $v_H=\sum_{x\in\mathcal{N}}v_H(x)\lambda_x\in V_H$,  $v:=v_H+\sum_{x\in\mathcal{N}}\left(v_H(x)-(\Inodal v_H)(x)\right)b_x\in V$ has the desired properties (for suitably chosen $b_x$).
The interpolation and support properties are obvious. The stability follows from
\begin{eqnarray*}
 \|\nabla v\|^2&\leq &C\left(\|\nabla v_H\|^2+\sum_{x\in\mathcal{N}}\left|v_H(x)-(\Inodal v_H)(x)\right|^2\|\nabla b_x\|^2\right)\\
&\leq &C\left(\|\nabla v_H\|^2+\Cinv'^2\sum_{x\in\mathcal{N}}\left|v_H(x)-(\Inodal v_H)(x)\right|^2\|\nabla \lambda_x\|^2\right)\\
&\leq &C\left(\|\nabla v_H\|^2+C'\Cinv'^2\sum_{T\in\triH}\|v_H-\Inodal v_H\|_{L^2(T)}^2 H_T^{-2}\right)\\
&\leq &C\left(\|\nabla v_H\|^2+C'\Cint^2\Cinv'^2\sum_{T\in\triH}\|\nabla v_H\|_{L^2(\omega_T)}^2\right)\\
&\leq & \Cinv^2 \|\nabla v_H\|^2,
\end{eqnarray*}
where we use $\|\nabla \lambda_x\|^2\approx |\support \lambda_x|^{(d-2)}$, the inverse inequality $\|v_H-\Inodal v_H\|_{L^\infty(T)}^2\lesssim H_T^{-d}\|v_H-\Inodal v_H\|_{L^2(T)}^2$, and \eqref{e:interr}.
\end{proof}
\medskip

In the forthcoming derivation of our method, the interpolation operator \eqref{e:clement} may be replaced by any linear bounded surjective operator that satisfies \eqref{e:interr}--\eqref{e:invert}. Hereby, \eqref{e:invert} may be relaxed in the sense that $\support v$ is not necessarily a subset of $\support v_H$ but that ${\support v\setminus\support v_H}$ covers at most a fixed (small) number of element layers about $\support v_H$.

\subsection{Multiscale Splitting and Modified Nodal Basis}\label{ss:multiscale}
Let $\Inodal: V\rightarrow V_H$ be a quasi interpolation operator according to the previous subsection. Then the kernel of $\Inodal$
\begin{equation*}\label{e:finescale}
 \Vf:=\{v\in V\;\vert\;\Inodal v=0\}
\end{equation*}
represents the microscopic features of $V$, i.e., all features that are not captured by $V_H$.
Given $v\in V_H$, define $\Tfine v\in\Vf$ by
\begin{equation*}\label{e:finescaleproj}
 a(\Tfine v,w)=a(v,w)\quad\text{for all }w\in \Vf.
\end{equation*}
The finescale projection operator $\Tfine:V_H\rightarrow \Vf$ leads to an orthogonal splitting with respect to the scalar product $a$
\begin{equation*}\label{e:splitting}
 V=\VmsH\oplus \Vf\quad\text{where}\quad \VmsH:=(V_H-\Tfine V_H).
\end{equation*}
Hence, any function $u\in V$ can be decomposed into $\ums\in \VmsH$ and $\ufine\in \Vf$, $u=\ums+\ufine$,
with $a(\ums,\ufine)=0$. Since $\dim\VmsH=\dim V_H$, the space $\VmsH$ can be regarded as a modified coarse space. The superscript ``ms'' abbreviates  ``multiscale'' and indicates that $\VmsH$, in addition, contains fine scale information.
The corresponding Galerkin approximation $\ums\in \VmsH$ satisfies
\begin{equation}\label{e:msmodel}
a( \ums,v) = G(v) \quad\text{ for all }v\in \VmsH.
\end{equation}
The error $(u-\ums)$ of the above method \eqref{e:msmodel} is analyzed in Section \ref{ss:errorms}.

Finally, we shall introduce a basis of $\VmsH$. The image of the nodal basis function $\lambda_x$ under the fine scale projection $\Tfine$ is denoted by $\phi_x=\Tfine\lambda_x\in \Vf$, i.e., $\phi_x$ satisfies the corrector problem
\begin{equation}\label{e:Tnodal}
 a(\phi_x,w)=a(\lambda_x,w)\quad\text{for all }w\in \Vf.
\end{equation}
We emphasize that the corrector problem is posed in the fine scale space $\Vf$, i.e., test and trial functions satisfy the constraint that their interpolation with respect to the coarse mesh vanishes.

A basis of $\VmsH$ is then given by the modified nodal basis
\begin{equation}\label{e:basiscoarse}
 \{\lambda_x-\phi_x\;\vert\; x\in\mathcal{N}\}.
\end{equation}
In general, the corrections $\phi_x$ of nodal basis functions $\lambda_x$, $x\in\mathcal{N}$, have global support, a fact which limits the practical use of the modified basis \eqref{e:basiscoarse} and the corresponding method \eqref{e:msmodel}.

\subsection{Localization}\label{ss:localization}
In Section~\ref{ss:errormslocal}, we will show that the correction $\phi_x$ decays exponentially fast away from $x$. Hence, simple truncation of the corrector problems to local patches of coarse elements yields localized basis functions with good approximation properties.

Let $k\in\mathbb{N}$. Define nodal patches of $k$-th order $\omega_{x,k}$ about $x\in\mathcal{N}$ by
\begin{equation}\label{e:omega}
 \begin{aligned}
 \omega_{x,1}&:=\support \lambda_x=\inter\left(\cup\left\{T\in\triH\;|\;x\in T\right\}\right),\\
 \omega_{x,k}&:=\inter\left(\cup\left\{T\in\triH\;|\;T\cap {\overline\omega}_{x,{k-1}}\neq\emptyset\right\}\right),\quad k=2,3,4\ldots .
\end{aligned}
\end{equation}
Define localized finescale spaces $\Vf(\omega_{x,k}):=\{v\in\Vf\;\vert\;v\vert_{\Omega\setminus\omega_{x,k}}=0\}$, $x\in\mathcal{N}$, by intersecting $\Vf$ with those functions that vanish outside the patch $\omega_{x,k}$.
The solutions $\phi_{x,k}\in \Vf(\omega_{x,k})$ of
\begin{equation}\label{e:Tnodallocal}
 a(\phi_{x,k},w)=a(\lambda_x,w)\quad\text{for all }w\in \Vf(\omega_{x,k}),
\end{equation}
are approximations of $\phi_x$ from \eqref{e:Tnodal} with local support.

We define localized multiscale finite element spaces
\begin{subequations}\label{e:modeldiscrete}
\end{subequations}
\begin{equation}
\VmsHk=\operatorname*{span}\{\lambda_x-\phi_{x,k}\;\vert\;x\in\mathcal{N}\}\subset V.
\tag{\ref{e:modeldiscrete}.a}
\end{equation}
The corresponding multiscale approximation of \eqref{e:model} reads: find $\umsloc{k}\in\VmsHk$ such that
\begin{equation}
a(\umsloc{k},v) = G(v) \quad\text{ for all }v\in\VmsHk.
\tag{\ref{e:modeldiscrete}.b}
\end{equation}
Note that $\dimension\VmsHk=|\mathcal{N}|=\dimension V_H$, i.e., the number of degrees of freedom of the proposed method \eqref{e:modeldiscrete} is the same as for the classical method \eqref{e:courantb}. The basis functions of the multiscale method have local support. The overlap is proportional to the parameter $k$. The error analysis of Section \ref{ss:errormslocal} suggests to choose $k\approx\log\tfrac{1}{H}$.

\begin{remark}\label{rem:furtherloc}
The localized modified basis functions could be localized further to vertex patches $\omega_x$, $x\in\mathcal{N}$, by simply multiplying them with the classical nodal basis functions; for any $x\in\mathcal{N}$ and any $y\in\mathcal{N}\cap\omega_{x,k}$, define $\phi_{x}^y:=\lambda_y \phi_{x,k}$. The generalized finite element space which is spanned by those $\mathcal{O}\left(  \left(\log(1/H)\right)^{d}\right)$ local basis functions per vertex has similar approximation properties as $V_{H,k}^{ms}$ (see \cite{Babuska96thepartition}).
\end{remark}

\section{Error Analysis}\label{s:error}
This section analyzes the proposed multiscale method in two steps. First, Subsection \ref{ss:errorms} presents an error bound for the idealized method \eqref{e:msmodel}. Then, Subsection \ref{ss:errormslocal} bounds the error of truncation to local patches and proves the main result, that is, an error bound for the multiscale method \eqref{e:modeldiscrete}.

As usual, the error analysis depends on the constant $\rho>0$ which represents shape regularity of the finite element mesh $\triH$;
\begin{equation}
\label{e:shapereg}
\rho:=\max_{T\in\triH}\rho_T\quad\text{with}\quad \rho_T:=\frac{\diam{B_T}}{\diam{T}}\; \text{ for } T\in\triH,
\end{equation}
where $B_T$ denotes the largest ball contained in $T$.

\subsection{Discretization Error}\label{ss:errorms}
\begin{lemma}\label{l:discerror}
Let $u\in V$ solve \eqref{e:model} and $\ums\in \VmsH$ solve \eqref{e:msmodel}. Then it holds
\begin{equation*}
\tnorm{u-\ums}\leq \Col^{1/2}\Cint\alpha^{-1/2}\|H g\|_{L^2(\Omega)}
\end{equation*}
with constants $\Col$ and $\Cint$ that only depend on $\rho$.
\end{lemma}
\begin{proof}
Recall the (local) approximation and stability properties \eqref{e:interr} of the interpolation operator $\Inodal$.
Due to the splitting from Section \ref{ss:multiscale}, it holds $u-\ums = \ufine$.
Since $\Inodal \ufine=0$, the application of \eqref{e:interr} and Young's inequality yield
\begin{align*}
 \tnorm{\ufine}^2& =G(\ufine) \leq \sum_{T\in\triH}\|g\|_{L^2(T)}\|\ufine-\Inodal \ufine\|_{L^2(T)}\\&\leq \frac{\Cint^2}{2\epsilon\alpha}\|H g\|_{L^2(\Omega)}^2+\frac{\epsilon}{2}\sum_{T\in\triH}\|A^{1/2}\nabla\ufine\|_{L^2(\omega_T)}^2
\end{align*}
for any $\epsilon>0$. Note that there exists a constant $\Col>0$ that only depends on $\rho$ such that the number of elements covered by  $\omega_T$ is uniformly bounded (w.r.t. $T$) by $\Col$. The choice $\epsilon=\Col^{-1}$ concludes the proof.
\end{proof}
\medskip

\begin{remark}
Substituting $\Inodal$ by the modified Cl\'ement interpolation operator presented in \cite{MR1736895} allows one to improve the error estimate in Lemma \eqref{l:discerror}. The term $\|H g\|_{L^2(\Omega)}$ can be replaced by data oscillations
$\left(\sum_{x\in\mathcal{N}}\|H (g-g_x)\|_{L^2(\omega_x)}^2\right)^{1/2}$ with some weighted averages $g_x$ of $g$ on $\omega_x$, $x\in\mathcal{N}$; we refer to \cite[Section 2]{MR1736895} for details. Further smoothness of the right hand side $g\in H^1(\Omega)$ then leads to quadratic convergence of the idealized method without localization.
\end{remark}

\subsection{Error of Localized Multiscale FEM}\label{ss:errormslocal}
First, we estimate the error due to truncation to local patches.
We will frequently make use of cut-off functions on element patches.
\begin{definition}\label{d:cutoff}
For $x\in\mathcal{N}$ and $m<M\in\mathbb{N}$, let $\eta^{m,M}_x :\Omega\rightarrow [0,1]$ be a continuous and weakly differentiable function such that
\begin{subequations}\label{e:cutoff}
\end{subequations}
\begin{align}
(\eta^{m,M}_x)\vert_{\omega_{x,m}}&=0,\tag{\ref{e:cutoff}.a}\label{e:cutoff.a}\\
(\eta^{m,M}_x)\vert_{\Omega\setminus \omega_{x,M}}&=1,\;\text{and}\tag{\ref{e:cutoff}.b}\label{e:cutoff.b}\\
\forall T\in\triH,\;\|\nabla\eta^{m,M}_x\|_{L^\infty(T)}&\leq \Cco(M-m)^{-1}H_T^{-1}\tag{\ref{e:cutoff}.c}\label{e:cutoff.c}
\end{align}
with some constant $\Cco$ that only depends on $\rho$. For example, one may choose $\eta_x^{m,M}\in S_H$ with nodal values
\begin{equation}\label{e:cutoffH}
\begin{aligned}
 \eta_x^{m,M}(x) &= 0\quad\text{for all }x\in\mathcal{N}\cap \omega_m,\\
 \eta_x^{m,M}(x) &= 1\quad\text{for all }x\in\mathcal{N}\cap \left(\Omega\setminus \omega_{x,M}\right),\text{ and}\\
 \eta_x^{m,M}(x)&= j(M-m)^{-1}\quad\text{for all }x\in\mathcal{N}\cap \partial\omega_{x,m+j},\;j=0,1,2,\ldots,M-m.
\end{aligned}
\end{equation}
\end{definition}

We prove the essential decay property of the corrector functions by some iterative Caccioppoli-type argument. Recall the notation $\tnorm{\cdot}_{\omega}:=\|A^{1/2}\nabla\cdot\|_{L^2(\omega)}$.
\begin{lemma}\label{l:basiserror}
For all $x\in\mathcal{N}$, $k,\ell \geq 2\in\mathbb{N}$, the estimate
\begin{equation*}
\tnorm{\phi_x-\phi_{x,\ell k}}\leq C_2\left(\frac{C_1}{\ell}\right)^{\frac{k-2}{2}}\tnorm{\phi_x}_{\omega_{x,\ell}}
\end{equation*}
holds with constants $C_1,C_2$ that only depend on $\rho$ and $\beta/\alpha$ but not on $x$, $k$, $\ell$, or $H$.
\end{lemma}
\begin{proof}
 Let $x\in\mathcal{N}$ and $\ell,k\geq 2\in\mathbb{N}$. Observe that
\begin{equation}\label{e:basiserror1}
\tnorm{\phi_x-\phi_{x,\ell k}}^2\leq \tnorm{\phi_x-v}^2=\tnorm{\phi_x-v}_{\omega_{x,\ell k}}^2+\tnorm{\phi_x}_{\Omega\setminus\omega_{x,\ell k}}^2,
\end{equation}
holds for all $v\in \Vf(\omega_{x,\ell k})$ using Galerkin orthogonality.

Let $\zeta_x:=1-\eta_x^{\ell(k-1)+1,\ell k-1}$ with
 a cutoff function $\eta_x^{\ell(k-1)+1,\ell k-1}$ as in Definition~\ref{d:cutoff}. According to \eqref{e:invert}, there exists $b_x\in V$ such that $\Inodal(b_x)=\Inodal(\zeta_x \phi_x)$, $\tnorm{b_x}\leq \Cinv \tnorm{\Inodal(\zeta_x \phi_x)}$, and $\support (b_x)\subset\omega_{x,\ell k}$. Hence, $v:=\zeta_x \phi_x-b_x\in\Vf(\omega_{x,\ell k})$  and
\begin{align*}
\tnorm{\phi_x-v}_{\omega_{x,\ell k}}&\leq \tnorm{\phi_x-\zeta_x \phi_x}_{\omega_{x,\ell k}\setminus \omega_{x,\ell (k-1)+1}}+\tnorm{b_x}_{\omega_{x,\ell k}\setminus \omega_{x,\ell (k-1)}}\\&\leq
\Cinv\Cint\left(\tnorm{\phi_x}_{\omega_{x,\ell k}\setminus \omega_{x,\ell (k-1)+1}}+\sqrt{\beta}\|\nabla(\zeta_x \phi_x)\|_{L^2(\omega_{x,\ell k}\setminus \omega_{x,\ell (k-1)})}\right).
\end{align*}
Since $\Inodal\phi_x=0$, the upper bound of the interpolation error \eqref{e:interr} and \eqref{e:cutoff.c} yield
\begin{align*}
 \|\nabla(\zeta_x \phi_x)\|&_{L^2(\omega_{x,\ell k}\setminus \omega_{x,\ell (k-1)})}^2
 \\
 &\leq C'''_2\hspace{-2ex}\sum_{T\in\triH:\;T\subset\overline{\omega}_{x,\ell k}\setminus \omega_{x,\ell (k-1)+1}}\hspace{-4ex}\left(H_T^2\|\grad \zeta_k\|^2_{L^\infty(T)}+\|\zeta_k\|^2_{L^\infty(T)}\right)\|\nabla \phi_x\|^2_{L^2(T)}\\&
\leq C''_2\alpha^{-1}\tnorm{\phi_x}_{\omega_{x,\ell k}\setminus \omega_{x,\ell (k-1)+1}}^2
\end{align*}
with $C''_2:=1+\Col \Cco^2\Cint^2$.
This leads to
\begin{equation}\label{e:basiserror2}
\tnorm{\phi_x-v}_{\omega_{x,\ell k}}\leq
C'_2\tnorm{\phi_x}_{\omega_{x,k \ell}\setminus \omega_{x,(k-1)\ell }},
\end{equation}
where $C'_2$ depends only on $\rho$ and $\sqrt{\beta/\alpha}$.
The combination of \eqref{e:basiserror1}, with $v=\zeta_x\phi_x-b_x$, and \eqref{e:basiserror2} yields
\begin{equation}\label{e:basiserror3}
\tnorm{\phi_x-\phi_{x,\ell k}}\leq C_2\tnorm{\phi_x}_{\Omega\setminus\omega_{x,\ell (k-1)}}.
\end{equation}

Further estimation of the right hand side in \eqref{e:basiserror3} is possible using cutoff functions $\eta_j:=\eta_x^{\ell(j-1)+1,\ell j}$ (cf. Definition~\ref{d:cutoff}), $j=2,3,\ldots,k-1$.
Observe that
\begin{equation}\label{e:est0}
\begin{aligned}
\|&A^{1/2}\grad \phi_x\|_{L^2(\Omega\setminus\omega_{x,\ell (k-1)})}^2\leq\|A^{1/2}\eta_{k-1}\grad \phi_x\|_{L^2(\Omega)}^2\\
&=\int_\Omega  (A\nabla\phi_x)\cdot\nabla(\eta_{k-1}^2
\phi_x)\;-\;2\int_\Omega\eta_{k-1}\phi_x (A\nabla\phi_x)\cdot\nabla
\eta_{k-1}.
\end{aligned}
\end{equation}
Let, according to Lemma~\ref{l:invert}, $b_{x,(k-1)}$ be chosen such that $\Inodal b_{x,(k-1)}=\Inodal(\eta_{k-1}^2
\phi_x)$. Then $\eta_{k-1}^2 \phi_x- b_{x,(k-1)}\in\Vf$.
Since $|\support(\grad\lambda_x)\cap\support(\eta_{k-1})|=0$ and $\support(\grad\eta_{k-1})=\omega_{x,(k-1)\ell}\setminus \omega_{x,(k-2)\ell+1}$, the first term on the right-hand side of \eqref{e:est0} can be rewritten as
\begin{equation}
 \begin{aligned}
\int_{\Omega}& (A\nabla\phi_x)\cdot\nabla(\eta_{k-1}^2
\phi_x)\\&
=\int_\Omega  (A\nabla\phi_x)\cdot \nabla(\eta_{k-1}^2
\phi_x- b_{x,(k-1)})+\int_\Omega  (A\nabla\phi_x)\cdot\nabla b_{x,(k-1)}\\
&=\int_\Omega  (A\nabla\phi_x)\cdot\nabla b_{x,(k-1)}\\
&\leq\Cinv\sqrt{\beta}\tnorm{\phi_x}_{\omega_{x,(k-1)\ell}\setminus \omega_{x,(k-2)\ell+1}}\|\nabla\Inodal(\eta_{k-1}^2
\phi_x)\|_{L^2(\omega_{x,(k-1)\ell}\setminus \omega_{x,(k-2)\ell+1})}.
\end{aligned}\label{e:est00}\end{equation}
With $\overline{\eta^2_T}:=|T|^{-1}\int_T\eta_{k-1}^2$ we have
\begin{align*}
\|&\nabla\Inodal(\eta_{k-1}^2
\phi_x)\|_{L^2(T)}=\|\nabla\Inodal((\eta_{k-1}^2-\overline{\eta^2_T})\phi_x)\|_{L^2(T)}\\
&\leq\Cint \|\nabla((\eta_{k-1}^2-\overline{\eta^2_T})\phi_x)\|_{L^2(T)}\\
&\leq\Cint \left(\|\eta_{k-1}^2-\overline{\eta^2_T}\|_{L^{\infty}(T)}\|\nabla\phi_x\|_{L^2(T)}+\|\nabla(\eta_{k-1}^2)\|_{L^{\infty}(T)}\|\phi_x\|_{L^2(T)}\right)\\
&\leq2\Cint\|\nabla(\eta_{k-1})\|_{L^{\infty}(T)} \left(\alpha^{-1/2}\diam(T)\tnorm{\phi_x}_{T}+\|\phi_x-\Inodal(\phi_x)\|_{L^2(T)}\right).
\end{align*}
Thus, the property \eqref{e:cutoff.c} of the cutoff function and the upper bound of the interpolation error \eqref{e:interr} yield
\begin{align}\label{e:est0a}
\tnorm{\Inodal(\eta_{k-1}^2
\phi_x)}_{\omega_{x,(k-1)\ell}\setminus \omega_{x,(k-2)\ell+1}}\leq C_1'\ell^{-1} \|A^{1/2}\grad\phi_x\|_{L^2(\Omega\setminus \omega_{x,(k-2)\ell})},
\end{align}
where $C_1'$ only depends on $\Cint$, $\Cco$, $\Col$, and $\sqrt{\beta/\alpha}$.
The same arguments allow one to bound the second term on the right-hand side in \eqref{e:est0}
\begin{equation}
 \begin{aligned}\label{e:est0b}
2\int_\Omega&\eta_{k-1}\phi_x (A\nabla\phi_x)\cdot\nabla
\eta_{k-1}\\&\leq2\hspace{-2ex}\sum_{T\in\triH:\;T\subset\overline{\omega}_{x,(k-1)\ell}\setminus \omega_{x,(k-2)\ell+1}}\hspace{-4ex}\|\grad \eta_{k-1}\|_{L^\infty(T)}\|A^{1/2} \grad \phi_x\|_{L^2(T)}\|A^{1/2}\phi_x\|_{L^2(T)}\\
&\leq C_1''\ell^{-1}\|A^{1/2}\grad\phi_x\|_{L^2(\Omega\setminus \omega_{x,(k-2)\ell})}^2,
\end{aligned}
\end{equation}
where $C_1''$ only depends on $\Cint$, $\Cco$, and $\sqrt{\beta/\alpha}$.
The combination of \eqref{e:est0}--\eqref{e:est0b} yields
\begin{equation}\label{e:est1}
\tnorm{\phi_x}_{\Omega\setminus\omega_{x,(k-1)\ell}}^2\leq C_1 \ell^{-1}\tnorm{\phi_x}_{\Omega\setminus \omega_{x,(k-2)\ell}}^2,
\end{equation}
where $C_1:=C_1'+C_1''$.
For $j=k-2,\ldots,2$, a similar argument (with $\eta_{k-1}$ replaced by $\eta_j$) yields
\begin{align}\label{e:estj}
\tnorm{\phi_x}_{\Omega\setminus \omega_{x,j\ell}}^2\leq C_1 \ell^{-1}\tnorm{\phi_x}_{\Omega\setminus \omega_{x,(j-1)\ell}}^2.
\end{align}
Starting from \eqref{e:est1}, the successive application of \eqref{e:estj} for $j=k-2,k-3,\ldots,2$ proves
\begin{equation}\label{e:decay}
\tnorm{\phi_x}_{\Omega\setminus\omega_{x,(k-1)\ell}}^2\leq (C_1 \ell^{-1})^{k-2}\tnorm{\phi_x}_{\omega_{x,\ell}}^2.
\end{equation}
Combining \eqref{e:basiserror3} and \eqref{e:decay}, we finally obtain the assertion.
\end{proof}
\medskip

\begin{lemma}\label{l:trace}
There is a constant $C_3$ that depends only on $\rho$ and $\beta/\alpha$, but not on $|\mathcal{N}|$, $k$, or $\ell$ such that
$$
\Tnorm{\sum_{x\in\mathcal{N}} v(x)(\phi_x-\phi_{x,\ell k})}^2
  \leq C_3 (\ell k)^d\sum_{x\in\mathcal{N}} v^2(x)\tnorm{\phi_x-\phi_{x,\ell k}}^2.
$$
\end{lemma}
\begin{proof}
For $x\in\mathcal{N}$, let $\zeta_x=1-\eta_x^{\ell k+1,\ell k+2}$ (cf. Definition~\ref{d:cutoff}). By Lemma \ref{l:invert} there exists a function $b_x\in V$ such that for any $w\in \Vf$ it holds $\Inodal b_x=\Inodal ((1-\zeta_x) w)$, $\text{supp}(b_x)\subset \text{supp}(\Inodal ((1-\zeta_x) w))\subset \omega_{x,\ell k+3}\setminus\omega_{x,\ell k}$, and $\tnorm{ b_x}_{\omega_{x,\ell k+3}\setminus\omega_{x,\ell k}}\leq \Cinv\tnorm{ \Inodal ((1-\zeta_x) w)}_{\omega_{x,\ell k+3}\setminus\omega_{x,\ell k}}$. We note that $w-\zeta_x w-b_x\in \Vf$ with support outside $\omega_{x,\ell k}$ , i.e.,~$a(\phi_x,w-\zeta_x w-b_x)=a(\lambda_x,w-\zeta_x w-b_x)=0$ and $a(\phi_{x,\ell k},w-\zeta_x w-b_x)=0$. With $w=\sum_{x\in\mathcal{N}}v(x)(\phi_x-\phi_{x,\ell k})\in \Vf$ we have
\begin{align*}
\tnorm{w}^2&=
\sum_{x\in\mathcal{N}} v(x)\, a(\phi_x-\phi_{x,\ell k},\zeta_x w+b_x)
\\
&\leq \sqrt{\beta}\sum_{x\in\mathcal{N}} |v(x)|\,\tnorm{\phi_x-\phi_{x,\ell k}} \cdot\|\nabla(\zeta_x w)\|_{L^2(\Omega)}\\
&\qquad +\sqrt{\beta}\sum_{x\in\mathcal{N}} |v(x)|\,\tnorm{\phi_x-\phi_{x,\ell k}} \cdot\Cinv\|\nabla(\Inodal((1-\zeta_x) w))\|_{L^2(\omega_{x,\ell k+3})}\\
&\leq 2\sqrt{\beta}\Cinv\Cint \sum_{x\in\mathcal{N}} |v(x)|\,\tnorm{\phi_x-\phi_{x,\ell k}} \cdot\|\nabla(\zeta_x w)\|_{L^2(\Omega)}\\
&\qquad +2\sqrt{\beta}\Cinv\Cint \sum_{x\in\mathcal{N}} |v(x)|\,\tnorm{\phi_x-\phi_{x,\ell k}} \cdot\|\nabla w\|_{L^2(\omega_{x,\ell k+4})}\\
&\leq 2\sqrt{\beta}\Cinv\Cint\sum_{x\in\mathcal{N}} |v(x)|\,\tnorm{\phi_x-\phi_{x,\ell k}} \cdot \,\|(\nabla\zeta_x)(1-\Inodal)w)\|_{L^2(\omega_{x,\ell k+2})}\\
&\quad +2\sqrt{\tfrac{\beta}{\alpha}}\Cinv\Cint\sum_{x\in\mathcal{N}} |v(x)|\,\tnorm{\phi_x-\phi_{x,\ell k}}\cdot\tnorm{w}_{\omega_{x,\ell k+4}}\\
&\leq 4\sqrt{\tfrac{\beta}{\alpha}}\Cinv\Cint^2\Cco\sum_{x\in\mathcal{N}} |v(x)|\,\tnorm{\phi_x-\phi_{x,\ell k}} \cdot\tnorm{w}_{\omega_{x,\ell k+4}}\\
&\leq 4\sqrt{\tfrac{\beta}{\alpha}}\Cinv\Cint^2\Cco\Cov(\ell k)^{d/2}\left(\sum_{x\in\mathcal{N}}v^2(x)\tnorm{\phi_x-\phi_{x,\ell k}}^2\right)^{1/2}\tnorm{w},
\end{align*}
where $\Cov (\ell k)^d$ represents an upper bound on the number of patches $\omega_{x,\ell k}$ that overlap a single element in the mesh. The result follows by dividing by $\tnorm{w}$ on both sides.
\end{proof}
\medskip

\begin{theorem}\label{thm3}
Let $u\in V$ solve \eqref{e:model} and, given $\ell,k\geq 2\in\mathbb{N}$, let $\umsloc{\ell k}\in V_{H,\ell k}^{ms}$ solve
\eqref{e:modeldiscrete}. Then
\begin{multline*}
 \tnorm{u-\umsloc{\ell k}}\leq C_4\|H_T^{-1}\|_{L^\infty(\Omega)}\left(\ell k\right)^{d/2}(C_1/\ell)^{\frac{k-2}{2}}\|g\|_{H^{-1}(\Omega)}\\+\Col^{1/2}\Cint\alpha^{-1/2}\|H g\|_{L^2(\Omega)}.
\end{multline*}
holds with $C_1$ from Lemma \ref{l:basiserror} and a constant $C_4$ that depends on $\alpha$, $\beta$ and $\rho$ but not on $H$, $k$, $\ell$, $g$, or $u$.
\end{theorem}
\begin{proof}
Let $\umsloctilde{\ell k}:=\sum_{x\in\mathcal{N}}\ums(x)\left(\lambda_x-\phi_{x,\ell k}\right)$, where $\ums(x)$, $x\in\mathcal{N}$, are the coefficients in the basis representation of $\ums$.
 Due to Galerkin orthogonality, Lemma \ref{l:discerror}, Lemma \ref{l:trace}, and the triangle inequality,
\begin{equation}\label{e:error1}
\begin{aligned}
 \tnorm{u-\umsloc{\ell k}}&\leq\tnorm{u-\umsloctilde{\ell k}}=\tnorm{u-\ums+\ums-\umsloctilde{\ell k}}\\ &\leq \Col^{1/2}\Cint\alpha^{-1/2}\|H g\|_{L^2(\Omega)}+\tnorm{\ums-\umsloctilde{\ell k}}.
\end{aligned}
\end{equation}
The application of Lemma \ref{l:basiserror} yields
\begin{align*}
\tnorm{\ums-\umsloctilde{\ell k}}^2&\leq C_3(\ell k)^{d}\sum_{x\in\mathcal{N}}\ums(x)^2\tnorm{\phi_x-\phi_{x,\ell k}}^2\\&\leq C_3(\ell k)^{d}C_2^2(C_1/\ell)^{k-2} \sum_{x\in\mathcal{N}}\ums(x)^2 \tnorm{\phi_x}_{\omega_{x,\ell}}^2.
\end{align*}
Furthermore, we have
\begin{align*}
\sum_{x\in\mathcal{N}}\ums(x)^2&\tnorm{\phi_x}^2_{\omega_{x,l}}
\leq \beta\Cinverse \sum_{T\in\mathcal{T}}H_T^{-2}\sum_{x\in T\cap\mathcal{N}}\ums(x)^2\|\lambda_x\|^2_{L^2(T)}\\
&\leq \beta\Cinverse'  \sum_{T\in\mathcal{T}}H_T^{-2}\biggl\|\sum_{x\in T\cap\mathcal{N}}\ums(x)\lambda_x\biggr\|^2_{L^2(T)}\\
&= \beta\Cinverse' \biggl\|H^{-2}\sum_{x\in\mathcal{N}}\ums(x)\lambda_x\biggr\|^2_{L^2(\Omega)}\\
&\leq \beta\Cinverse'\biggl(\|H^{-2}\ums\|^2_{L^2(\Omega)}+ \biggl\|H^{-2}\sum_{x\in\mathcal{N}}\ums(x)(\phi_x-\Inodal\phi_x)\biggr\|^2_{L^2(\Omega)}\biggr)\\
&\leq \tfrac{\beta}{\alpha}\Cinverse'(C_{\text{F}}\|H_T^{-2}\|_{L^\infty(\Omega)}+\Cint)\tnorm{\ums}^2.
\end{align*}
where $\Cinverse$ and $\Cinverse'$ depend on $\rho$ and $C_{\text{F}}=C_{\text{F}}(\Omega)$ is the constant from Friedrichs' inequality. This yields
\begin{equation}\label{e:error3}
\begin{aligned}
\tnorm{\ums-\umsloctilde{\ell k}}&\leq C_4' \|H_T^{-1}\|_{L^\infty(\Omega)}(\ell k)^{d/2}(C_1/\ell)^{(k-2)/2} \tnorm{\ums}\\&\leq C_4\|H_T^{-1}\|_{L^\infty(\Omega)}(\ell k)^{d/2}(C_1/\ell)^{(k-2)/2} \|g\|_{H^{-1}(\Omega)},
\end{aligned}
\end{equation}
where $C_4$ only depends on $C_2$, $C_3$, $C'_{\text{inv}}$, $C_{\text{F}}$, $\Cint$, and $\sqrt{\beta}/\alpha$. The assertion follows readily by combining \eqref{e:error1} and \eqref{e:error3}.
\end{proof}
\begin{remark}
The error estimate in Theorem~\ref{thm3} contains a factor $\|H^{-1}\|_{L^\infty(\Omega)}$. However, its influence on the total error can be controlled by choosing the localization parameter $k$ proportional to $\log(1/\|H^{-1}\|_{L^\infty(\Omega)})$. For non-uniform meshes, it is recommended to vary the choice of the localization parameter in space according to $k\approx\log\tfrac{1}{H}$. We neglect this opportunity to avoid overloading the paper.
\end{remark}

\medskip

\section{Discretization of the Fine Scale Computations}\label{s:disc}
In this section, we focus on how to compute numerical approximations to the local basis functions $\lambda_x-\phi_{x,\ell k}$ and thereby to the multiscale solution $\umsloc{\ell k}$. In order to do this, we need to extend the error analysis of Section 3 to a fully discrete setting. There is a lot of freedom in choosing different finite elements and different refinement strategies, see e.g.~\cite{MR2319044,MR2553176}. We will focus on a very simple and natural approach. We assume that the local basis functions are computed using subgrids of a fine scale reference mesh, which is a (possibly space adaptive) refinement of the coarse grid $\triH$.

More precisely, let $\trih$ be the result of one uniform refinement and several conforming but possibly non-uniform refinements of the coarse mesh $\triH$. We introduce $h:\overline\Omega\rightarrow \mathbb{R}_{>0}$ as the $\trih$-piecewise constant mesh width function with $h_t:=h|_t=\text{diam}(t)$ for all $t\in\trih$. We construct the finite element space
\begin{equation*}\label{e:courantfineglob}
S_{h}:=\{v\in C^0(\Omega)\;\vert\;\forall t\in\trih(\Omega),v\vert_t \text{ is a polynomial of total degree}\leq 1\}.
\end{equation*}
We let $u_h\in V_h:=S_h\cap H^1_0(\Omega)$ be the reference solution that satisfies
\begin{equation}\label{e:eiscref}
a(u_h,v)=G(v) \quad\text{for all }v\in V_h.
\end{equation}

Locally on each patch we let
\begin{equation}
\Vhxf(\omega_{x,k}):=\Vf(\omega_{x,k})\cap V_h=\{v\in V_{h}\;\vert\; \Inodal v=0\text{ and } v\vert_{\Omega\setminus\omega_{x,k}}=0\}.
\end{equation}
The numerical approximation $\phi_{x,k}^h\in \Vhxf(\omega_{x,k})$ of the corrector $\phi^h_{x,k}$ is determined by
\begin{equation*}
a(\phi_{x,k}^h,w)=a(\lambda_x,w) \quad\text{for all }w\in \Vhxf(\omega_{x,k}).
\end{equation*}
We denote the discrete multiscale finite element space
\begin{equation*}
\VmsHkh=\operatorname*{span}\{\lambda_x-\phi^h_{x,k}\;\vert\;x\in\mathcal{N}\}.
\end{equation*}
The corresponding discrete multiscale approximation $\umsloch{k}\in \VmsHkh$ fulfills
\begin{equation}\label{e:modelnested}
a(\umsloch{k},v) = G(v) \quad\text{ for all }v\in \VmsHkh.
\end{equation}

\begin{theorem}\label{thm5}
 Let $u\in V$ solve \eqref{e:model} and let $\umsloch{\ell k}\in \VmsHkh$  solve
\eqref{e:modelnested}. Then
\begin{multline*}
\tnorm{u-\umsloch{\ell k}}\leq \tilde{C}_4\|H_T^{-1}\|_{L^\infty(\Omega)}\left(\ell k\right)^{d/2}(\tilde{C}_1/\ell)^{\frac{k-2}{2}}\|g\|_{H^{-1}(\Omega)}\\+\Col^{1/2}\Cint\alpha^{-1/2}\|H g\|_{L^2(\Omega)} + \tnorm{ u-u_h},
\end{multline*}
where $\tilde{C}_4$ only depends on $rho$, $\alpha$ and $\beta$.
\end{theorem}

\begin{remark}[Multiscale splitting by nodal interpolation]
Having discretized the fine scale computation, i.e., having replaced the infinite dimensional space $V$ by some finite element space $V_h\subset C^0(\overline{\Omega})$ we are allowed to replace the Cl\'ement-type interpolation by classical nodal interpolation. This leads to the variational multiscale method in \cite{MRX}, which is a modification of the method first presented in \cite{MR2553176}. Because nodal interpolation satisfies the conditions \eqref{e:interr}--\eqref{e:invert}, Theorem~\ref{thm5} establishes an a priori error bound for the multiscale method \cite{MRX}. However, the constant $\Cint$ in \eqref{e:interr} depends on the ratio $H/h$ of the discretization scales if $d>1$ ($\Cint\approx \log(H/h)$ in $2d$ and $\Cint\approx (H/h)^{-1}$ in $3d$, c.f. \cite{MR853662}). Hence, for nodal interpolation, the constants $\tilde{C}_1,\tilde{C}_4$ in Theorem~\ref{thm5} depend on $H/h$ in a similar fashion. In 2d this can still be acceptable because the dependence on $H/h$ is only logarithmic.
\end{remark}

\begin{remark}[Estimates for the fine scale error]
The finite element space $V_h$ may be replaced by any finite element space that contains $V_h$, e.g., by piecewise polynomials of higher order.
The third part in the error bound in Theorem \ref{thm5} can be bounded in terms of data, mesh parameter $h$, and polynomial degree using standard a priori error estimates. For example, if $A\in W^{1,\infty}(\Omega)$ (bounded with bounded weak derivative) and $\varepsilon$ is the smallest present scale, i.e., $\|\nabla A\|_{L^{\infty}(\Omega)}\lesssim \varepsilon^{-1}$, the third term in the error bound in Theorem \ref{thm5} may be replaced by the worst case bound $C h \varepsilon^{-1}$ for a first-order ansatz space $V_h$ (see \cite{PS11}). It is shown in \cite{PS11} that for highly varying but smooth coefficient $A$, higher order ansatz spaces are superior.
\end{remark}

\begin{remark}[Periodic coefficient]
Let $\Omega$ be some square or cube, $g\in L^2(\Omega)$, $A$ be smooth and periodic, $A(x)=A(x/\varepsilon)$, with some small scale parameter $\varepsilon>0$, and let $u_\varepsilon$ denote the corresponding solution of \eqref{e:model}. Choose uniform meshes $\triH$ and $\trih$ with $H>\varepsilon>h$ and $k\approx\log(H^{-1})$. With regard to the previous comment, Theorem~\ref{thm5} yields the error bound
\begin{equation*}
 \tnorm{u_\varepsilon-\umsloch{\ell k}}\leq C_g(H+\tfrac{h}{\varepsilon}).
\end{equation*}
With $h\sim \varepsilon H$ the error in the approximation becomes independent of the fine scale oscillations without any so-called resonance effects as they are observed, e.g., in \cite{MR1455261}. We emphasize that periodicity can be exploited to reduce the number  corrector problems to be solved essentially.
\end{remark}

\begin{remark}[Solution of the local problems]
The local problems need to be solved in the spaces $\Vhxf(\omega_{x,k})$. This is a standard finite element space with the additional constraint that the trial and test functions should have no component in $V_H$. In practice this constraint is realized using Lagrange multipliers.

The resulting coarse scale system of equations is of the same size as the original problem, $\text{dim}(\VmsHkh)=\text{dim}V_H$ and it is still sparse. The number of non-zero entries will be larger and depend on $k$. Note however that the non-zero entries in the stiffness matrix decay exponentially away from the diagonal.
\end{remark}

\begin{proof}[Proof of Theorem~\ref{thm5}]
We use the triangle inequality $\tnorm{u-u^{\text{ms},h}_{H,\ell k}}\leq \tnorm{u-u_h}+\tnorm{u_h-u^{\text{ms},h}_{H}}+\tnorm{u^{\text{ms},h}_{H}-u^{\text{ms},h}_{H,\ell k}}$ and follow the arguments from the proof of Theorem~\ref{thm3} simply replacing $V$ by $V_h$ and using  Lemmas \ref{l:discerrorh}, \ref{l:basiserrorh}, and \ref{l:traceh} below (discrete versions of Lemmas \ref{l:discerror}, \ref{l:basiserror}, and \ref{l:trace}) to bound the last two terms.
\end{proof}

\begin{lemma}[Discrete version of Lemma~\ref{l:discerror}]\label{l:discerrorh}
Let $u_h\in V_h$ solve \eqref{e:eiscref} and $\umsh\in \VmsHh$ solve \eqref{e:modelnested} with $k$ large enough so that $\omega_{x,k}=\Omega$ for all $x\in\mathcal{N}$. Then
\begin{equation*}
\tnorm{u_h-\umsh}\leq \Col^{1/2}\Cint\alpha^{-1/2}\|H g\|_{L^2(\Omega)}
\end{equation*}
holds with constants $\Col$ and $\Cint$ that only depend on $\rho$.
\end{lemma}
\begin{proof}
Note that $\ufineh:=u_h-\umsh$ is the unique element of $\Vf_h:=\Vf\cap V_h$ such that $a(\ufineh,v)=G(v)$ for all $v\in \Vf_h$. The Lemma follows from the same arguments in the proof of Lemma \ref{l:discerror}.
\end{proof}
\medskip

In the remaining part of this Section, $A\lesssim B$ abbreviates an inequality $A\leq C\,B$ with some generic constant $0\leq C<\infty$ that does not depend on the mesh sizes $H$, $h$ and the localization parameters. The constant may depend on the contrast $\beta/\alpha$ but not on the geometrical or topological structure of the coefficient $A$.

To establish discrete versions of Lemmas~\ref{l:basiserror} and \ref{l:trace} we are facing the technical difficulty that the product of $v\in V_h$ and some cut-off function $\eta$ from Definition~\ref{d:cutoff} is not necessarily an element of $V_h$. However, the subsequent lemma shows that the product $\eta v$ can be approximated sufficiently well by elements from $V_h$.

\begin{lemma}\label{l:cutoffh}
For all $x\in\mathcal{N}$, $M>m\in\mathbb{N}$, and corresponding cut-off function $\eta^{m,M}_x$  defined in \eqref{e:cutoffH} there exists some $v\in \Vf_h(\omega_{x,M+1})$  such that
\begin{equation*}
\tnorm{\eta^{m,M}_x\phi^h_x-v}\lesssim  \frac{1}{M-m} \tnorm{\phi^h_x}_{\omega_{x,M+1}\setminus\omega_{x,m-1}}.
\end{equation*}
Furthermore the statement also holds if $\eta^{m,M}_x$ is replaced by $1-\eta^{m,M}_x$ and $v\in \Vf_h(\Omega\setminus\omega_{x,m-1})$.
\end{lemma}
\begin{proof}
Let $x\in\mathcal{N}$, $M>m\in\mathbb{N}$ be fixed and define $\eta:=\eta^{m,M}_x$.
Let $\Inodalh:V\cap C(\bar{\Omega})\rightarrow V^h$ be the nodal interpolant with respect to the mesh $\trih$. Recall its (local) approximation and stability properties
\begin{equation*}
 \|\nabla (v-\Inodalh v)\|_{L^2(t)}\lesssim h_t\|\nabla^2 v\|_{L^2(t)}\quad\text{and}\quad \|\Inodalh v\|_{L^2(t)}\lesssim \| v\|_{L^2(t)}
\end{equation*}
for all polynomials $v$. According to Lemma~\ref{l:invert}, there exists some $b_x\in V^h$ such that $\Inodal(b_x)=\Inodal(\Inodalh(\eta \phi^h_x))$, $\tnorm{b_x}\lesssim \tnorm{\Inodal(\Inodalh(\eta \phi^h_x))}$, and $\support (b_x)\subset\omega_{x,M+1}\setminus\omega_{x,m-1}$. Hence, $v:=\Inodalh(\eta \phi^h_x)-b_x\in\Vf_h(\omega_{x,M+1})$. Since $\Inodal \Inodalh\bar{\eta}_T\phi_x^h=\bar{\eta}_T\Inodal \phi_x^h=0$ for $\bar{\eta}_T=|T|^{-1}\int_T \eta$, we get
\begin{align*}
&\tnorm{\eta\phi^h_x-v}^2= \tnorm{\eta\phi^h_x-\Inodalh(\eta \phi^h_x)+b_x}^2 \\
&\lesssim \hspace{-0.5cm}\sum_{t\in\trih:t\subset \bar{\omega}_{x,M}\setminus\omega_{x,m}}\hspace{-0.45cm}\|\nabla(\eta\phi^h_x-\Inodalh(\eta \phi^h_x))\|^2_{L^2(t)}+\tnorm{\Inodal(\Inodalh((\eta-\bar{\eta}_T) \phi^h_x))}^2 \\
&\lesssim \hspace{-0.5cm}\sum_{t\in\trih:t\subset \bar{\omega}_{x,M}\setminus\omega_{x,m}}\hspace{-0.45cm}h_t^2\|\nabla^2(\eta\phi^h_x)\|^2_{L^2(t)}+\hspace{-0.7cm}\sum_{T\in\triH:T\subset \bar{\omega}_{x,M+1}\setminus\omega_{x,m-1}}\hspace{-0.5cm}H^{-2}_T\|\Inodalh((\eta-\bar{\eta}_T ) \phi^h_x)\|^2_{L^2(T)} \\
&\lesssim \hspace{-0.5cm}\sum_{t\in\trih:t\subset \bar{\omega}_{x,M}\setminus\omega_{x,m}}\hspace{-0.45cm}h_t^2\left(\|\nabla^2 \eta\|^2_{L^\infty(t)}\|\phi^h_x\|^2_{L^2(t)}+\|\nabla \eta\|^2_{L^\infty(t)}\|\nabla\phi^h_x\|^2_{L^2(t)}\right)\\&\qquad+\hspace{-0.7cm}\sum_{T\in\triH:T\subset \bar{\omega}_{x,M+1}\setminus\omega_{x,m-1}}\hspace{-0.5cm}H^{-2}_T \|\eta-\bar{\eta}_T\|^2_{L^{\infty}(T)}\|\phi^h_x\|^2_{L^2(T)}\\
&\lesssim (M-m)^{-1}\tnorm{\phi^h_x}_{\omega_{x,M+1}\setminus\omega_{x,m-1}}^2
\end{align*}
using the property \eqref{e:cutoff.c} of $\eta$ and Poincar\'e's inequality.
This proves the first part of the Lemma.

The second part concerning $1-\eta$ follows using the same argument but with $v\in \Vf_h(\Omega\setminus\omega_{x,m-1})$.
\end{proof}
\medskip

\begin{lemma}[Discrete version of Lemma~\ref{l:basiserror}]\label{l:basiserrorh}
For all $x\in\mathcal{N}$, $k,\ell \geq 2\in\mathbb{N}$ the estimate
\begin{equation*}
\tnorm{\phi^h_x-\phi^h_{x,\ell k}}\leq \tilde{C}_2\left(\frac{\tilde{C}_1}{\ell}\right)^{\frac{k-2}{2}}\tnorm{\phi^h_x}_{\omega_{x,\ell}}
\end{equation*}
holds with constants $\tilde{C}_1,\tilde{C}_2$ that only depend on $\rho$ and $\beta/\alpha$ but not on $x$, $k$, $\ell$, $h$, or $H$.
\end{lemma}
\begin{proof}
Let $\zeta_x:=1-\eta_x^{\ell(k-1)+1,\ell k-1}$ with $\eta_x^{\ell(k-1)+1,\ell k-1}$ as in equation \eqref{e:cutoffH} in Definition~\ref{d:cutoff}. Then there exists a $v\in \Vf_h(\omega_{x,\ell k})$ such that,
\begin{align*}
\tnorm{\phi^h_x-v}_{\omega_{x,\ell k}}&\leq \tnorm{\phi^h_x-\zeta_x \phi^h_x}_{\omega_{x,\ell k}}+ \tnorm{\zeta_x\phi^h_x-v}_{\omega_{x,\ell k}}\\
&\lesssim\tnorm{\phi^h_x}_{\omega_{x,\ell k}\setminus \omega_{x,\ell (k-1)+1}}+\tnorm{\zeta_x \phi^h_x}_{\omega_{x,\ell k-1}\setminus \omega_{x,\ell (k-1)+1}}.
\end{align*}
Furthermore, using the same argument as in Lemma \ref{l:basiserror},
\begin{equation*}\tnorm{\zeta_x \phi^h_x}_{\omega_{x,\ell k-1}\setminus \omega_{x,\ell (k-1)+1}}\lesssim \tnorm{\phi^h_x}_{\omega_{x,\ell k}\setminus \omega_{x,\ell (k-1)+1}}
\end{equation*}
which yields
\begin{equation*}
 \tnorm{\phi^h_x-\phi^h_{x,\ell k}}\leq \tnorm{\phi^h_x-v}_{\omega_{x,\ell k}}^2+\tnorm{\phi^h_x}_{\Omega\setminus\omega_{x,\ell k}}\lesssim \tnorm{\phi^h_x}_{\Omega\setminus \omega_{x,\ell (k-1)+1}}.
\end{equation*}

Now let $\eta_j:=\eta_x^{\ell(j-1)+1,\ell j}$ (cf. Definition~\ref{d:cutoff}), $j=2,3,\ldots,k-1$ and note that
\begin{equation*}
\|A^{1/2}\grad \phi^h_x\|_{L^2(\Omega\setminus\omega_{x,\ell (k-1)})}^2\leq a(\phi^h_x,\eta_{k-1}^2
\phi^h_x)-\;2\int_\Omega\eta_{k-1}\phi^h_x(A\nabla\phi^h_x)\cdot\nabla\eta_{k-1}.
\end{equation*}
The second term can be treated exactly as in Lemma \ref{l:basiserror} and, hence, bounded by $\ell^{-2}\tnorm{ \phi^h_x}_{\Omega\setminus\omega_{x,\ell (k-2)}}^2$. We make use of Lemma \ref{l:cutoffh} to bound the first term. There exists $v\in \Vf_h(\Omega\setminus\omega_{x,\ell (k-1)+1})$ such that
\begin{equation*}
a(\phi^h_x,\eta_{k-1}^2
\phi^h_x)\leq \tnorm{\phi_x^h}_{\omega_{\ell (k-1)}\setminus\omega_{\ell (k-2)+1}}\tnorm{\eta^2_{k-1}
\phi^h_x-v}\lesssim \ell^{-2}\tnorm{ \phi^h_x}_{\Omega\setminus\omega_{x,\ell (k-2)}}^2.
\end{equation*}
The final assertion follows by similar arguments as in the proof of Lemma~\ref{l:basiserror}.
\end{proof}
\medskip

\begin{lemma}[Discrete version of Lemma~\ref{l:trace}]\label{l:traceh}
There is a constant $\tilde{C}_3$ depending only on $\rho$ and $\beta/\alpha$, but not on $|\mathcal{N}|$, $k$, or $\ell$ such that,
\begin{equation*}
\biggl\|\biggl|\sum_{x\in\mathcal{N}} v(x)(\phi^h_x-\phi^h_{x,\ell k})\biggr|\biggr\|^2\leq \tilde{C}_3 (\ell k)^d\sum_{x\in\mathcal{N}} v^2(x)\tnorm{\phi^h_x-\phi^h_{x,\ell k}}^2.
\end{equation*}
\end{lemma}
\begin{proof}
For $x\in\mathcal{N}$, let $\zeta_x=1-\eta_x^{\ell k+1,\ell k+2}$ (cf. Definition~\ref{d:cutoff}) and $z=\sum_{x\in\mathcal{N}} v(x)\,(\phi_x-\phi_{x,\ell k})$. We have,
\begin{equation*}
\biggl\|\biggl|\sum_{x\in\mathcal{N}} v(x)\,(\phi_x-\phi_{x,\ell k})\biggr|\biggr\|^2=
\sum_{x\in\mathcal{N}} v(x)\, a(\phi^h_x-\phi^h_{x,\ell k},\zeta_x z+(1-\zeta_x) z)=\text{I}+\text{II}.
\end{equation*}
The first term $\text{I}:=\sum_{x\in\mathcal{N}} v(x)\, a(\phi^h_x-\phi^h_{x,\ell k},\zeta_x z)$ can be treated in exactly the same way as in the proof of Lemma \ref{l:trace}. We focus on the second term. Due to Lemma~\ref{l:cutoffh} there exists a $w\in \Vf_h(\Omega\setminus\omega_{\ell k})$ such that
\begin{align*}
\text{II}&:= \sum_{x\in\mathcal{N}} v(x)\, a(\phi^h_x-\phi^h_{x,\ell k},(1-\zeta_x) z-w)\\
&\lesssim \left(\sum_{x\in\mathcal{N}} |v(x)|^2\, \tnorm{\phi^h_x-\phi^h_{x,\ell k}}^2\right)^{1/2}\left(\sum_{x\in\mathcal{N}}\tnorm{(1-\zeta_x) z-w}^2\right)^{1/2}\\
&\lesssim \left(\sum_{x\in\mathcal{N}} |v(x)|^2\, \tnorm{\phi^h_x-\phi^h_{x,\ell k}}^2\right)^{1/2}\left(\sum_{x\in\mathcal{N}}\tnorm{z}_{\omega_{x,\ell k+2}\setminus\omega_{\ell k+1}}^2\right)^{1/2}\\
&\lesssim (\ell k)^{d/2}\left(\sum_{x\in\mathcal{N}} |v(x)|^2\, \tnorm{\phi^h_x-\phi^h_{x,\ell k}}^2\right)^{1/2}\tnorm{z}.
\end{align*}
The result follows immediately.
\end{proof}

\section{Numerical Experiments}\label{s:numexp}
Numerical experiments shall validate our theoretical results from the previous sections.
\subsection{Experimental setup}
\begin{figure}
\begin{center}
\includegraphics[height=0.225\textwidth]{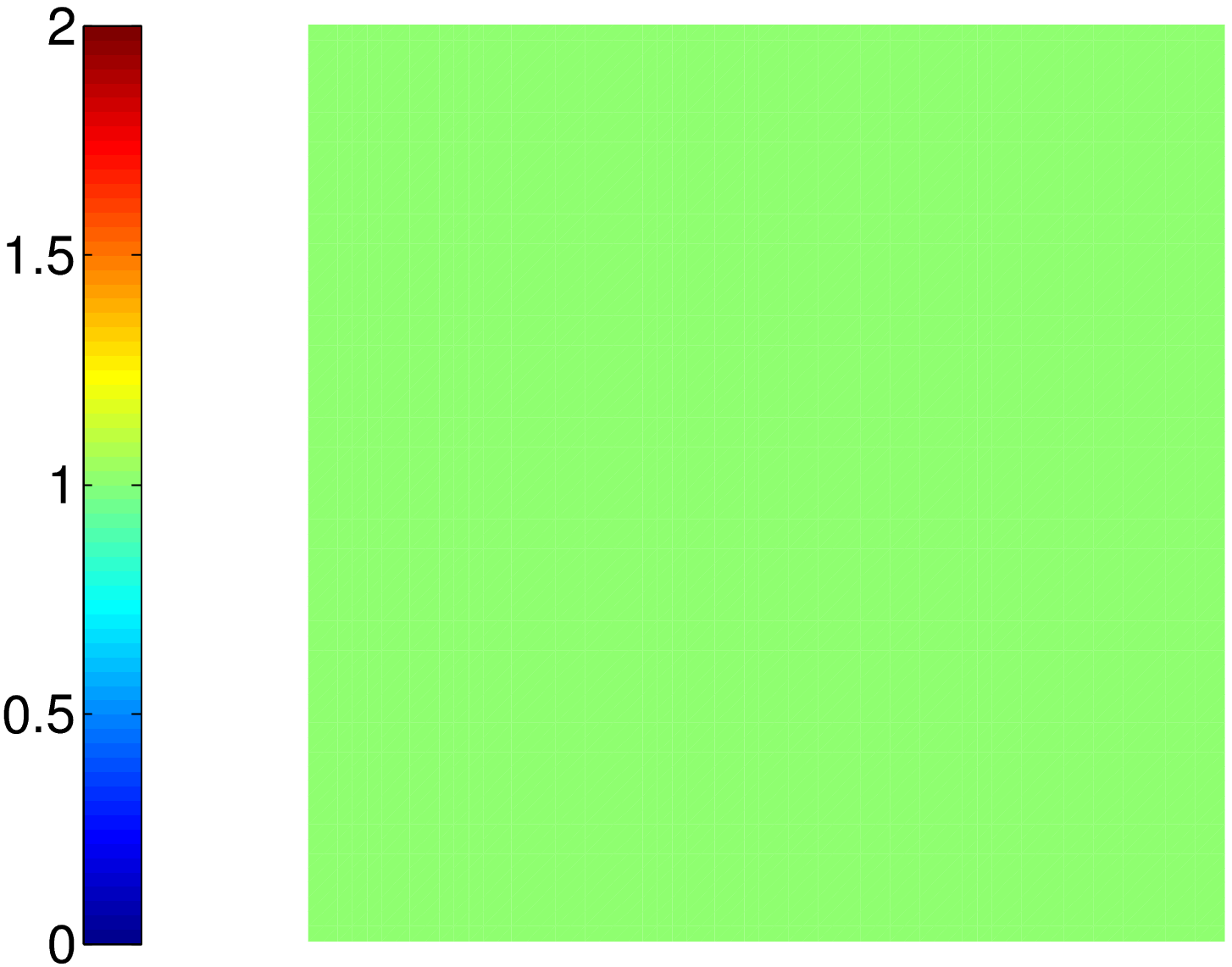}\hspace{0.03\textwidth}\includegraphics[height=0.22\textwidth]{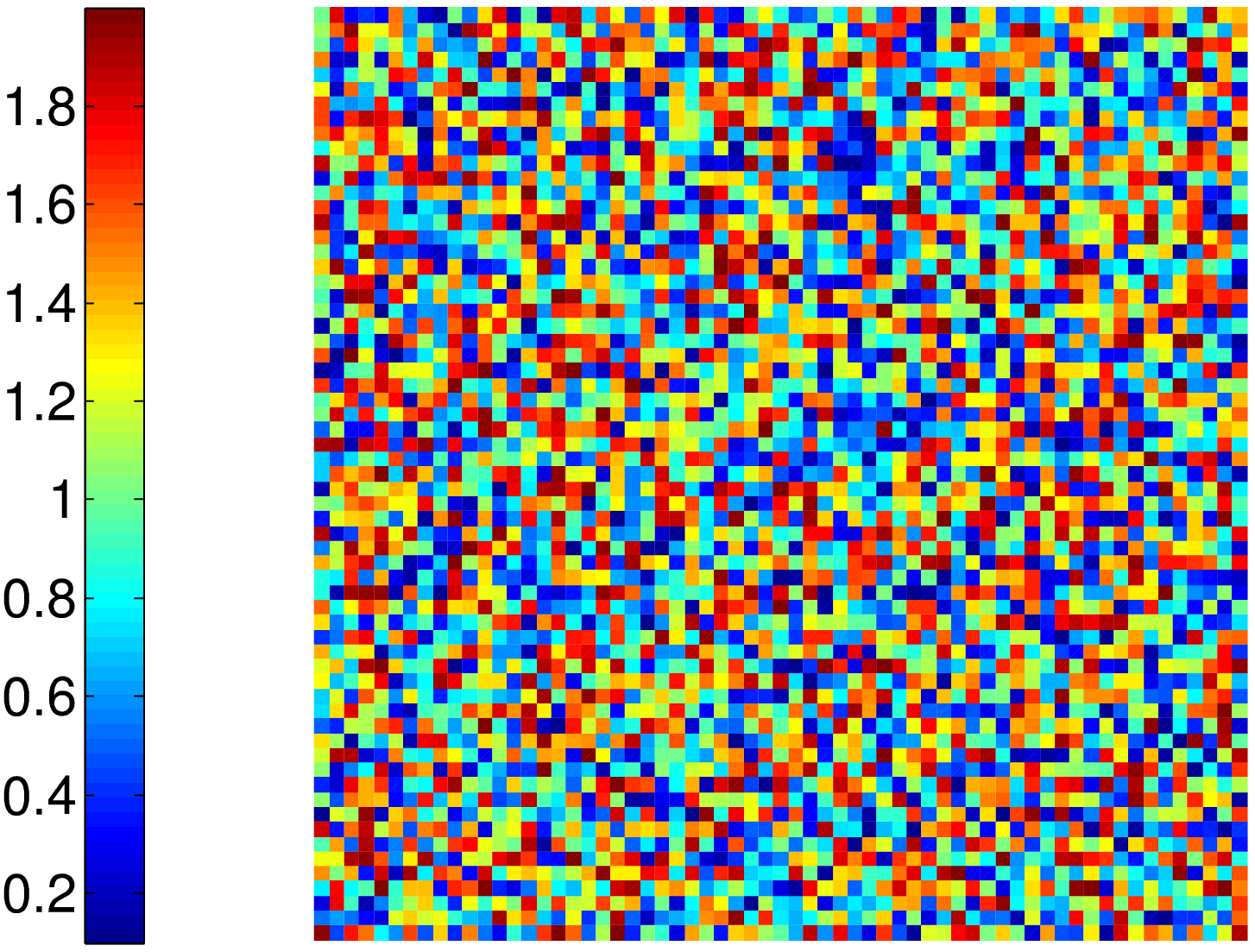}\hspace{0.03\textwidth}\includegraphics[height=0.22\textwidth]{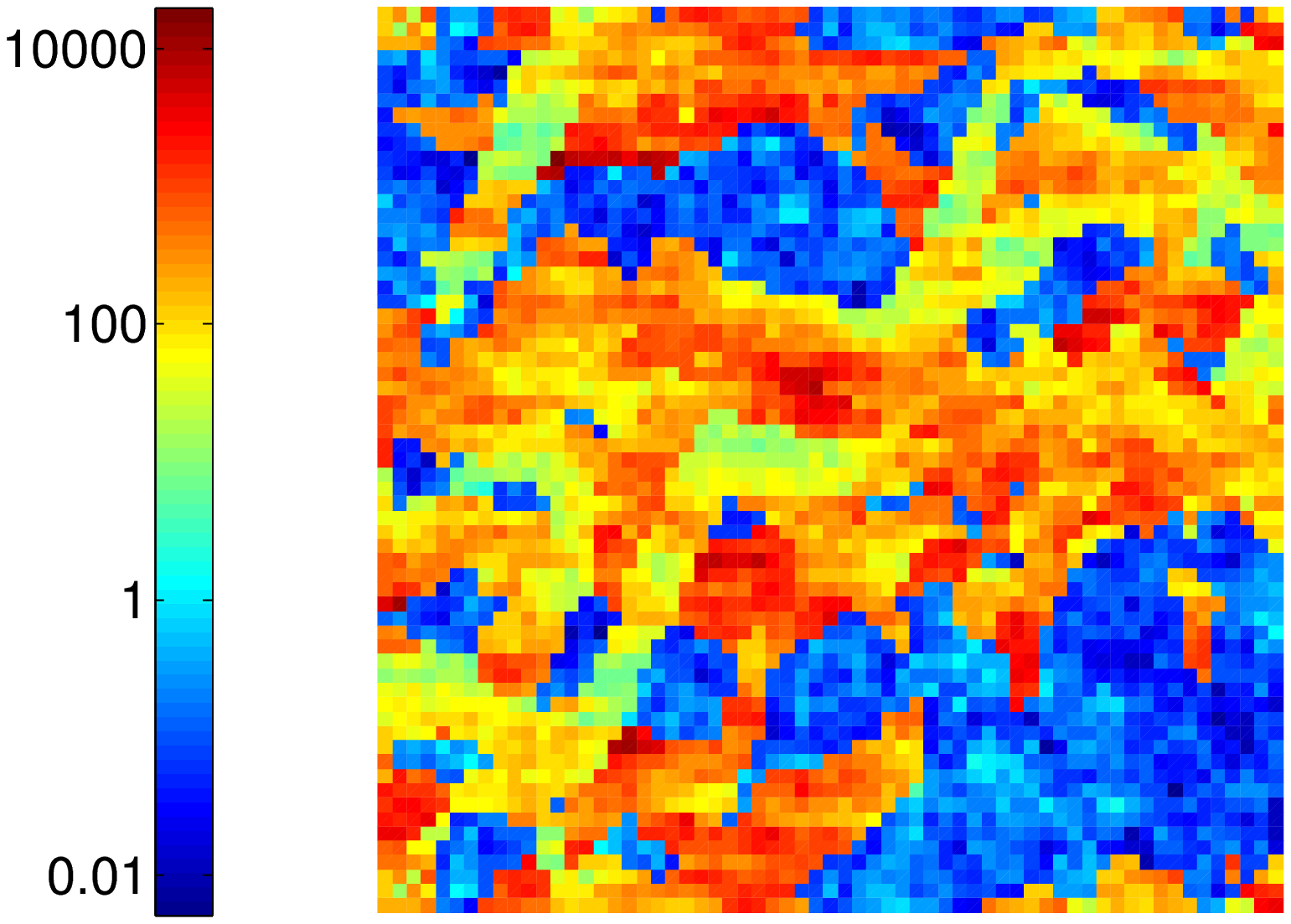}\end{center}
\caption{Scalar coefficient used in the numerical experiment: $A_1$ (left), $A_2$ (middle), $A_3$ (right).}\label{fig:coeff}
\end{figure}

\begin{figure}
\begin{center}
\includegraphics[width=0.25\textwidth]{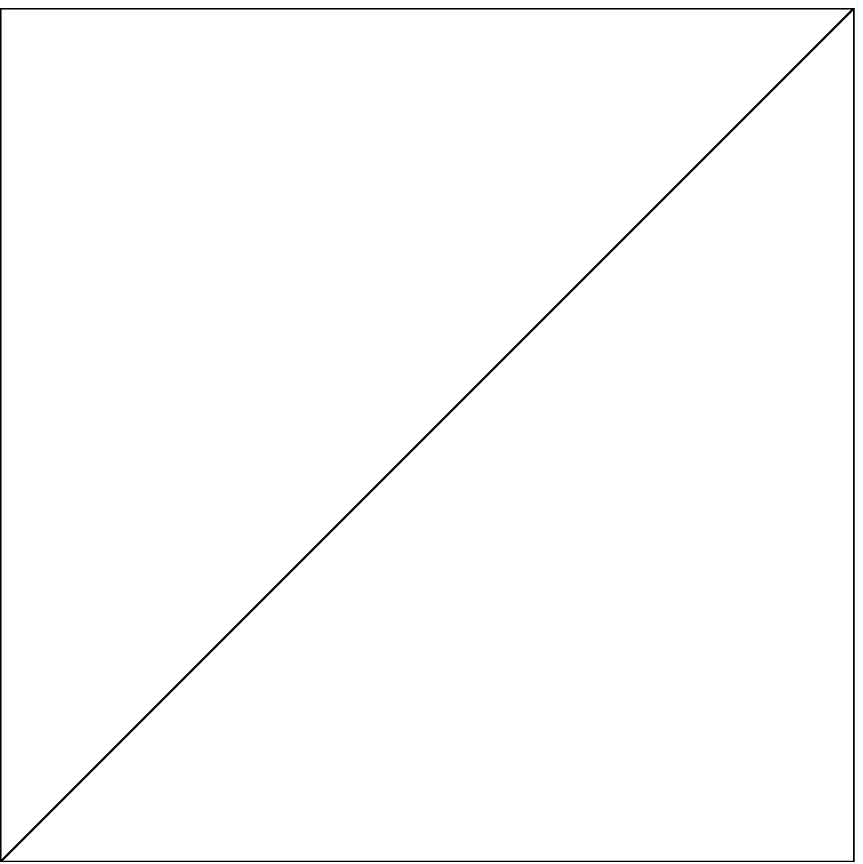}\hspace{2ex}\includegraphics[width=0.25\textwidth]{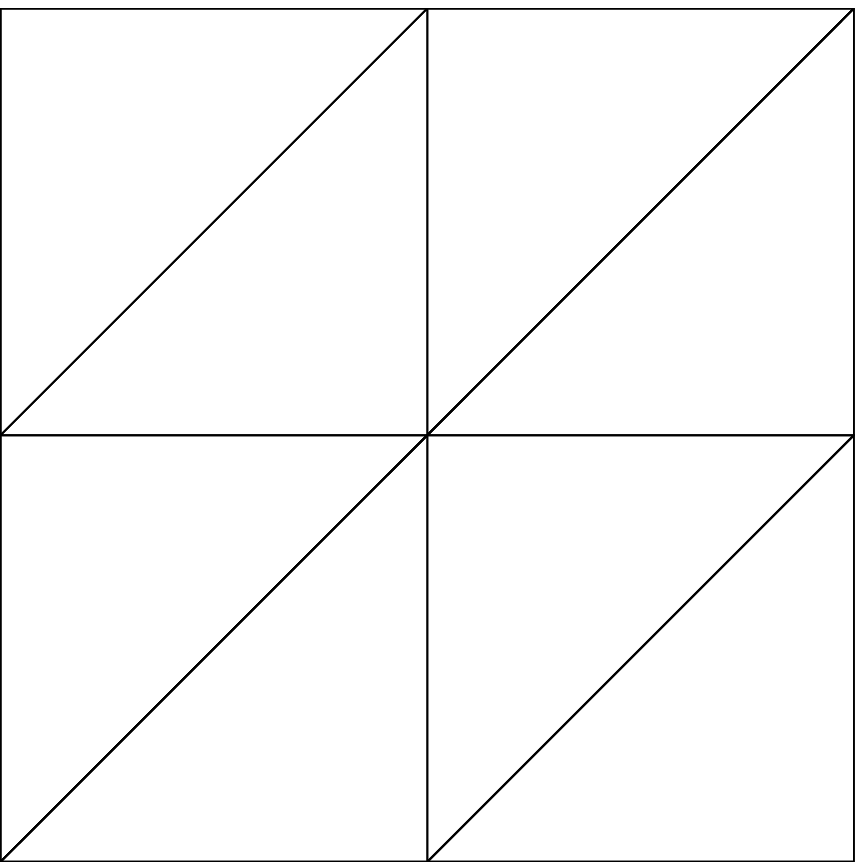}\hspace{2ex}\includegraphics[width=0.25\textwidth]{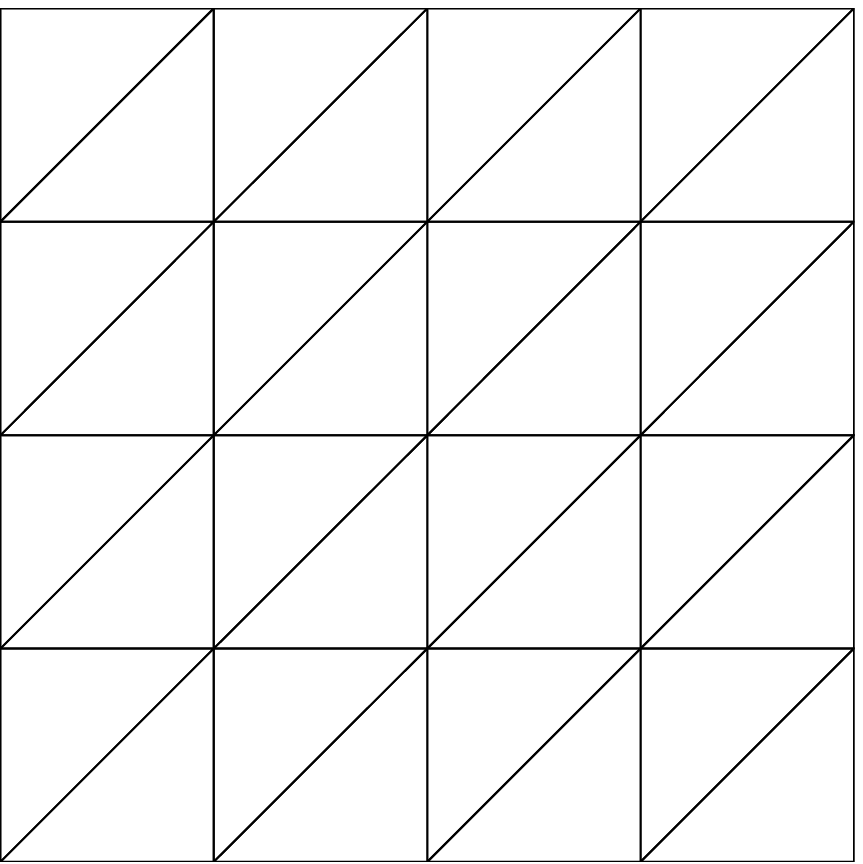}\end{center}
\caption{Uniform triangulations of the unit square.}\label{fig:meshes}
\end{figure}

Let $\Omega$ be the unit square and the outer force $g\equiv 1$ in $\Omega$. Consider three different choices for the scalar coefficient $A_1, A_2, A_3$ with increasing difficulty as depicted in Figure~\ref{fig:coeff}. The coefficient $A_1=1$ is constant. The coefficient $A_2$ is piecewise constant with respect to a uniform Cartesian grid of width $2^{-6}$. The values in each grid cell are chosen in the range $[1/20,2]$; the contrast $\beta(A_2)/\alpha(A_2)\leq 40$ is moderate. The coefficient $A_3$ is piecewise constant with respect to the same uniform Cartesian grid of width $2^{-6}$. Its values are taken from the data of the SPE10 benchmark, see \texttt{http://www.spe.org/web/csp/}. The contrast for $A_3$ is large, $\beta(A_3)/\alpha(A_3)\approx 4\cdot 10^6$.
Consider uniform coarse meshes of size $H=2^{-1},2^{-2},\ldots,2^{-6}$ of $\Omega$ as depicted in Figure~\ref{fig:meshes}. Note that none of these meshes resolves the rough coefficients $A_2$ and $A_3$ appropriately.

The reference mesh $\trih$ has width $h=2^{-9}$. Since no analytical solutions are available, the standard finite element approximation $u_h\in V_h$ on the reference mesh $\trih$ serves as the reference solution. All fine scale computations are performed on subsets of $\trih$.

The approximations are compared with this reference solution only. Doing this, we assume that $u_h$ is sufficiently accurate. True errors would behave similar in the beginning but level off at some point when the reference error $\tnorm{u-u_h}$ dominates the upscaling error.

\subsection{Results for the energy error}
\begin{figure}
\begin{center}
\includegraphics[height=0.25\textheight]{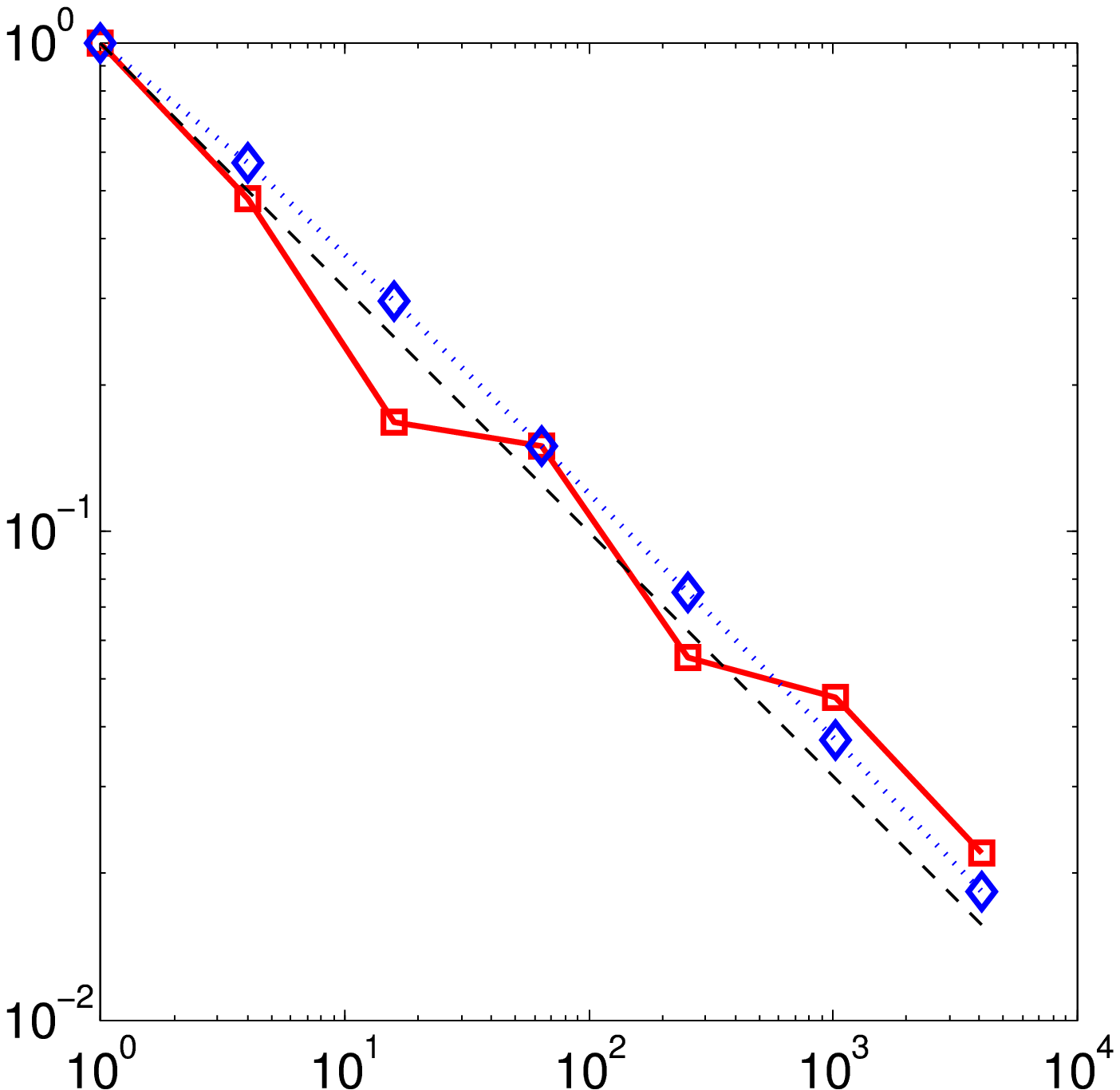}\vspace{1.5ex}\\
\includegraphics[height=0.25\textheight]{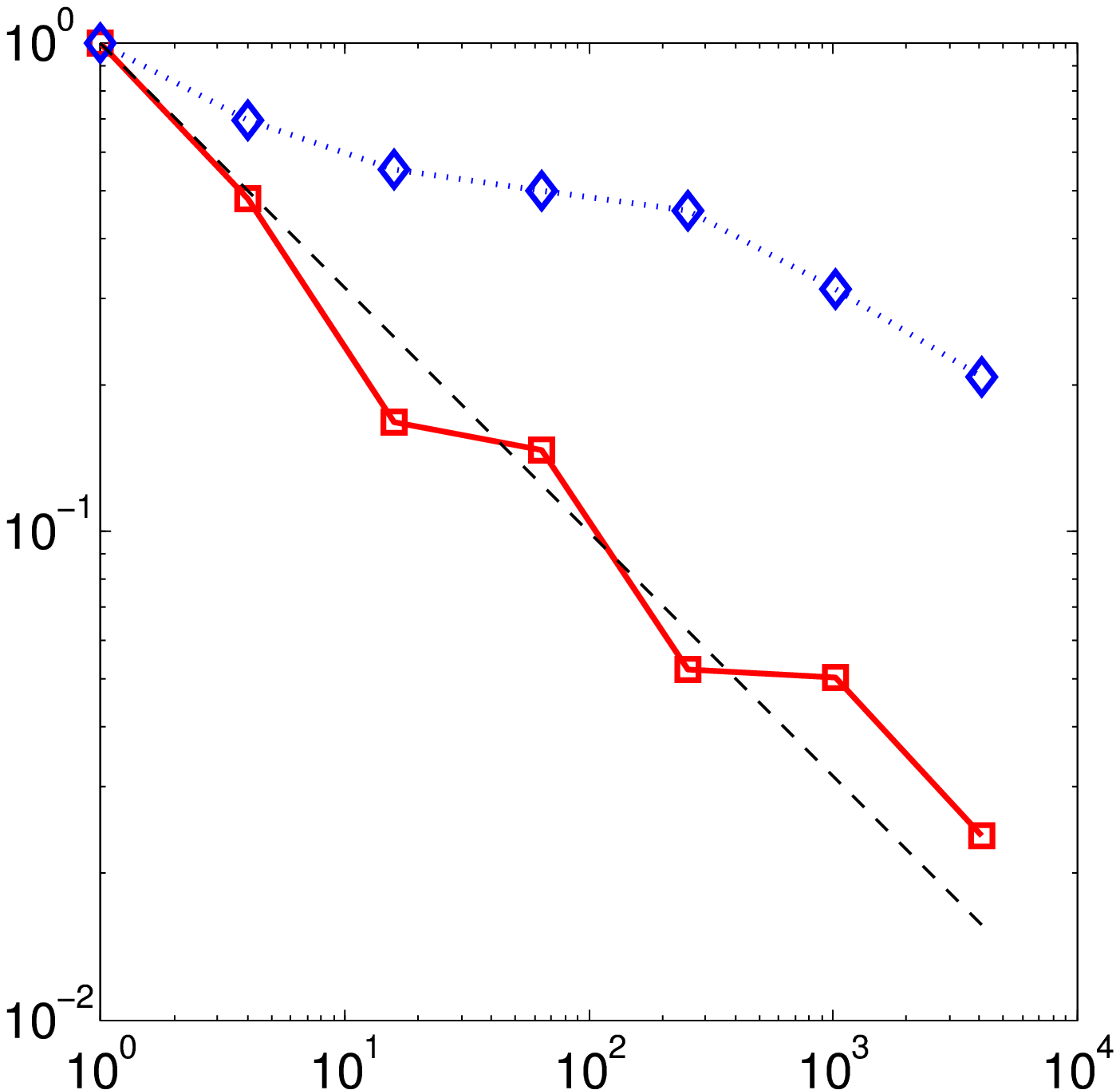}\vspace{1.5ex}\\
\includegraphics[height=0.25\textheight]{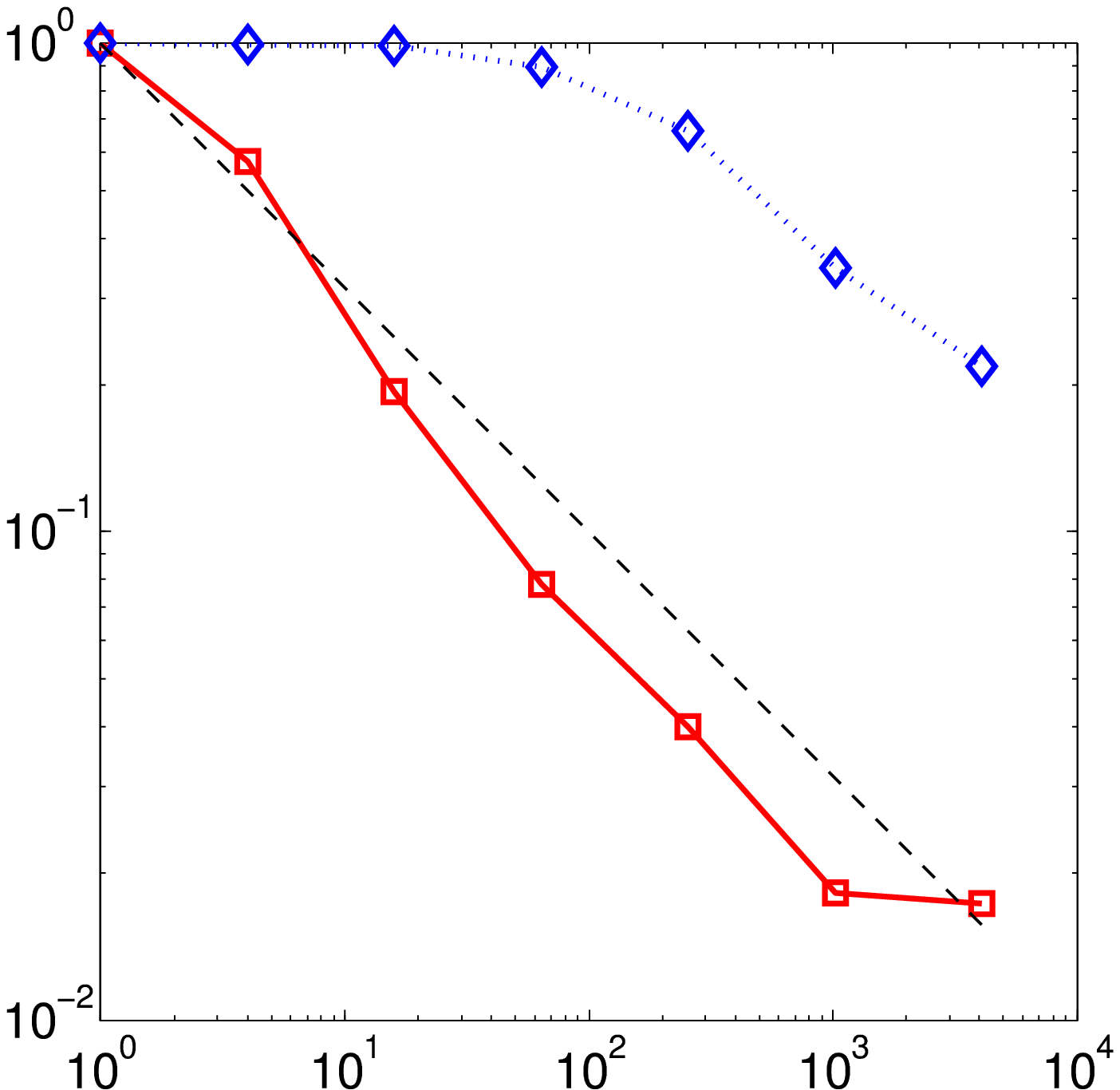}\end{center}
\caption{Relative energy errors $\tnorm{ u_h-\umsloch{k}}/\tnorm{ u_h}$ ($\square$ solid red) with localization parameter $k=\lceil 2\log(1/H)\rceil$ and $\tnorm{ u_h-u_H}/\tnorm{ u_h}$ ($\Diamond$ dotted blue) vs. number of degrees of freedom $N_{\operatorname*{dof}}\approx H^{-2}$ for different coefficients: $A_1$ (top), $A_2$ (middle), $A_3$ (bottom). The dashed black line is $N_{\operatorname*{dof}}^{-1/2}$.}\label{fig:H1}
\end{figure}

Figure~\ref{fig:H1} depicts the energy errors of the new multiscale method and the classical P1FEM (see \eqref{e:courantb}) with respect to the same coarse mesh. Depending on the coarse discretization scale $H$, the localization parameter $k$ is chosen to be $\lceil 2\log(1/H)\rceil$. The logarithmic dependence on $1/H$ is motivated by our a priori analysis. The choice of the constant $2$ is based on numerical tests. It turns out that, in all experiments, this choice leads to the desired linear textbook convergence (rate $-1/2$) of the energy error (w.r.t. to the number of degrees of freedom $N_{\operatorname*{dof}}=|\mathcal{N}|\approx H^{-2}$) related to the sequence of multiscale approximations. Pre-asymptotic effects are not observed.
In particular, the performance of our method does not seem to be affected by the high contrast present in $A_3$. Whether our estimates on the decay of the corrector functions are sub-optimal or have worst-case character with respect to contrast is an issue of present research.

Observe that the classical P1FEM suffers from the lacks of approximability and regularity and converges only poorly for the rough coefficients $A_2$ and $A_3$.

\subsection{Results for the $L^2$ error}
\begin{figure}
\begin{center}
\includegraphics[height=0.25\textheight]{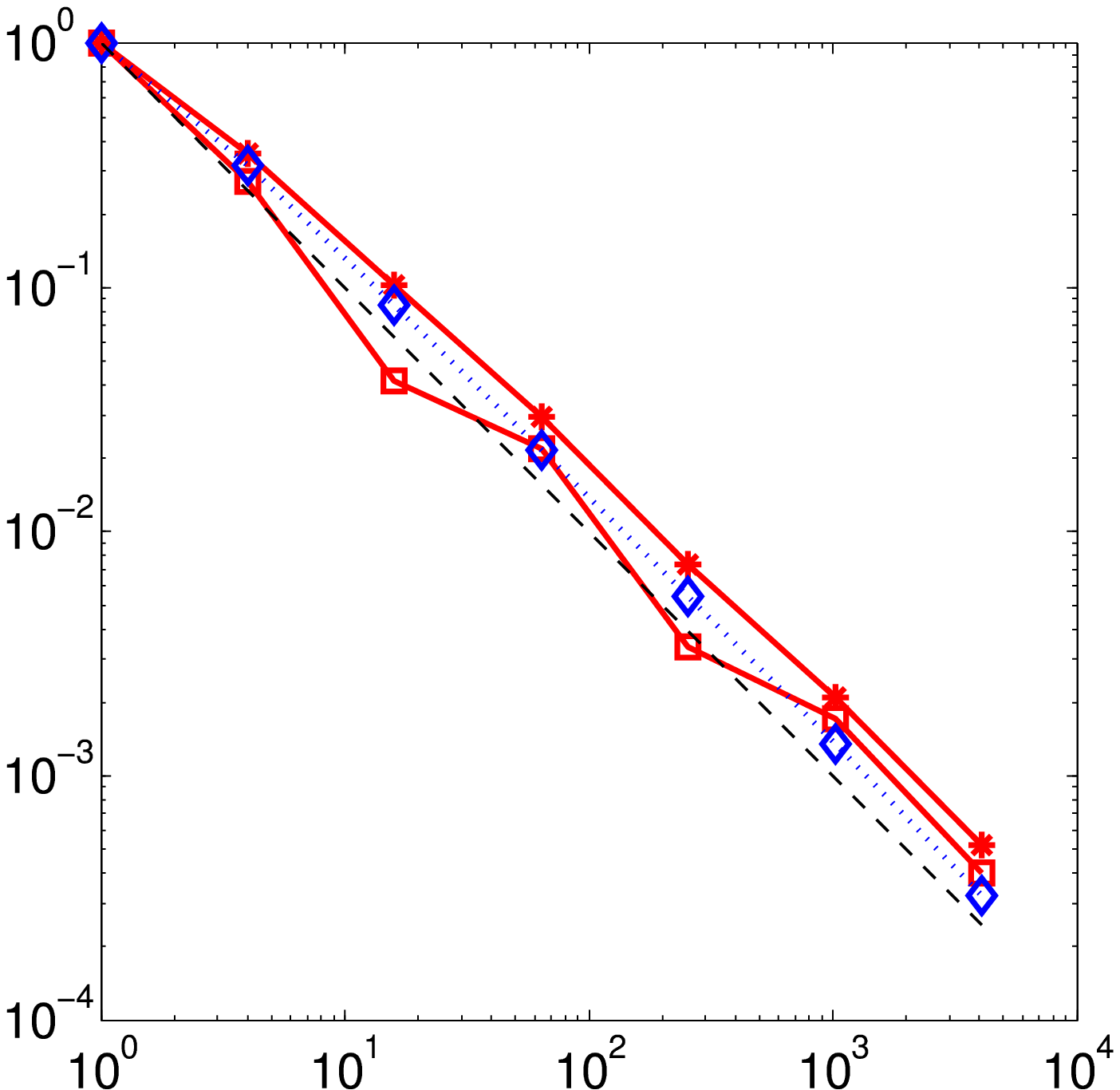}\vspace{1.5ex}\\
\includegraphics[height=0.25\textheight]{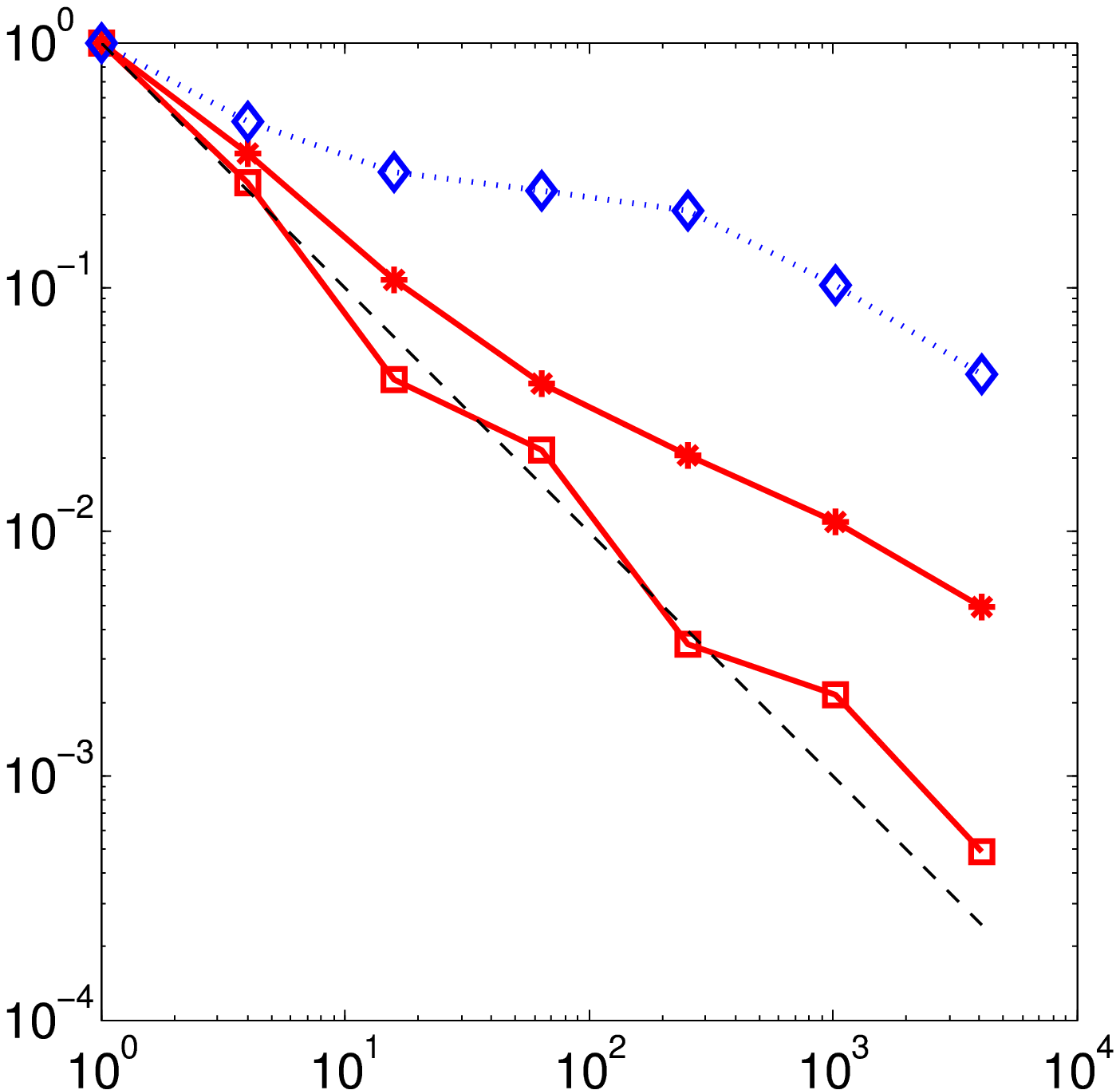}\vspace{1.5ex}\\
\includegraphics[height=0.25\textheight]{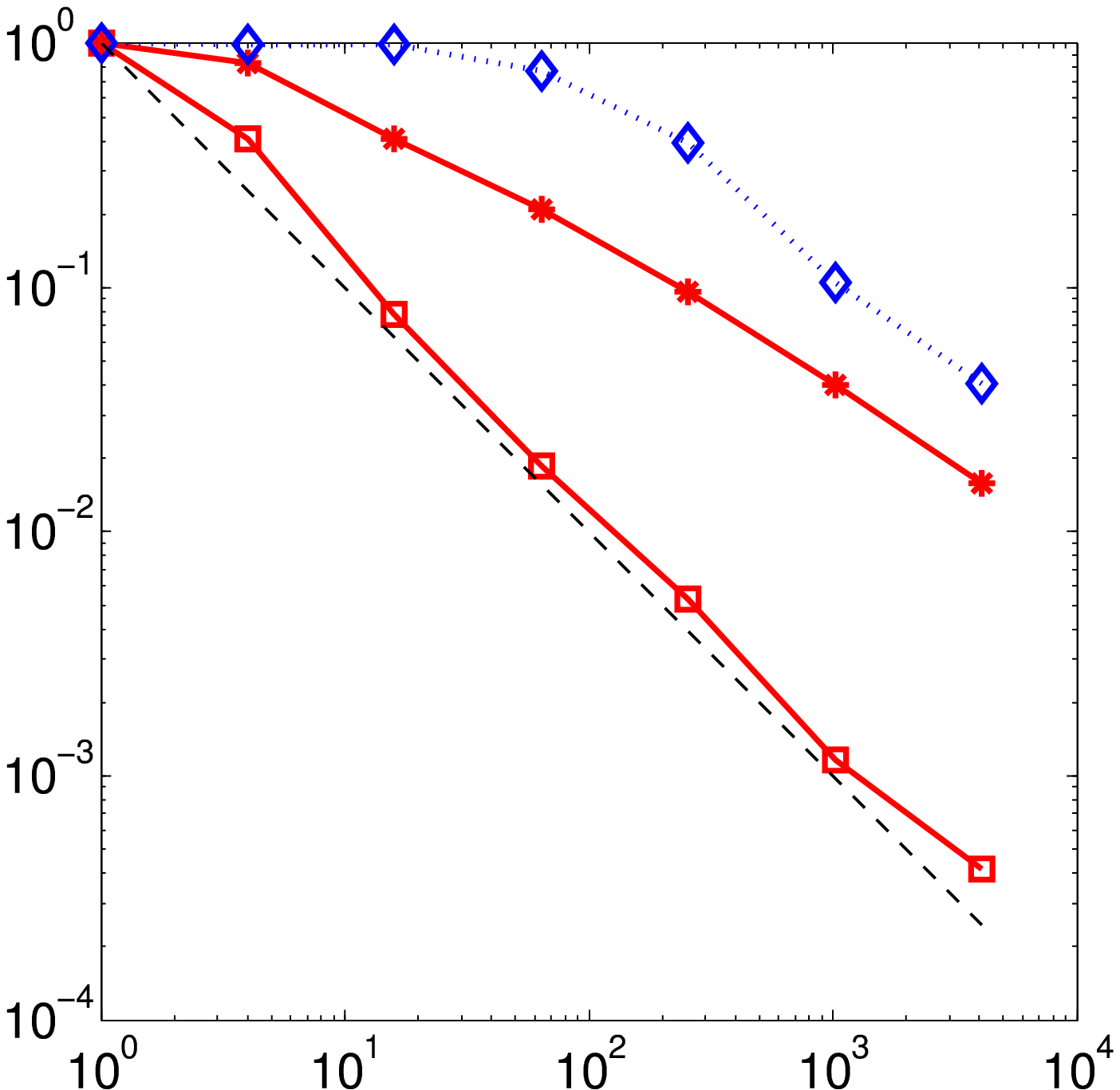}\end{center}
\caption{Relative $L^2$ errors $\| u_h-\umsloch{k}\|/\| u_h\|$ ($\square$ solid red), $\| u_h-\Inodal\umsloch{k}\|/\| u_h\|$ ($*$ solid red) with localization parameter $k=\lceil 2\log(1/H)\rceil$ and $\| u_h-u_H\|/\| u_h\|$ ($\Diamond$ dotted blue) vs. number of degrees of freedom $N_{\operatorname*{dof}}\approx H^{-2}$ for different coefficients: $A_1$ (top), $A_2$ (middle), $A_3$ (bottom). The dashed black line is $N_{\operatorname*{dof}}^{-1}$.}\label{fig:L2}
\end{figure}

Figure~\ref{fig:L2} shows $L^2$ errors of the new multiscale method and the classical P1FEM. Again, the choice of the localization parameter $k=\lceil 2\log(1/H)\rceil$ yields the optimal convergence rate $-1$ for our method in all experiments (w.r.t. to the number of degrees of freedom $N_{dof}=|\mathcal{N}|\approx H^{-2}$) without any pre-asymptotic behavior. This observation is justified by a standard Aubin-Nitsche duality argument. Define $e:=u_h-\umsloch{k}\in L^2(\Omega)$ and let $z_h\in V_h$ solve $$a(z_h,v_h)=\int_\Omega e v_h\quad\text{for all }v_h\in V_h.$$
Galerkin orthogonality leads to
$$\|u_h-\umsloch{k}\|_{L^2(\Omega)}^2=a(z_e-\zmsloch{k},e)\leq\tnorm{z_h-\zmsloch{k}}\tnorm{u_h-\umsloch{k}},$$
where $\zmsloch{k}\in\VmsHkh$ is the Galerkin projection of $z_h$ onto the discrete multiscale finite element space $\VmsHkh$.
The estimates for the energy error (see Section~\ref{s:disc}) and the present choice of $k$ yield the $L^2$ estimate
$$\|u_h-\umsloch{k}\|_{L^2(\Omega)}\lesssim H^2\|g\|_{L^2(\Omega)}.$$

More importantly, we observe that the $L^2$ error between $u_h$ and $\Inodal \umsloch{k}$ converges nicely at a rate close to $-3/4$ without pre-asymptotic effects. This is remarkable because, $\Inodal \umsloch{k}$ is a truly coarse approximation. $\Inodal \umsloch{k}$ is an element of the coarse $P1$ finite element space. Hence, it cannot capture microscopic features of the solution. The rate of convergence (with respect to the number of degrees of freedom) is limited by $\tfrac{1+s}{2}$ for some $s\in [0,1]$ which is related to the regularity of the solution ($u\in H^{1+s}$ for some $s\in [0,1]$). However, $\Inodal \umsloch{k}$ approximates the macroscopic behavior of the solution accurately with only very few degrees of freedom. Note that the storage complexity of the modified basis is of order $\mathcal{O}(h^{-2}\log{1/H})$ whereas its interpolation can be stored in $\mathcal{O}(H^{-2}\log{1/H})$. Once the coarse system matrix of the multiscale method is assembled, $\Inodal \umsloch{k}$ can be computed
without using any fine scale information from the modified basis whereas this would be required to represent the full multiscale approximation $\umsloch{k}$.

\section{Application to Multiscale Methods}\label{s:application}
In this section we discuss three multiscale methods and how the presented analysis relates to each of them.

\subsection{The Variational Multiscale Method}
The variational multiscale method was first introduced in \cite{MR1660141}. The function space $V$ is here split into a coarse part (standard finite element space on a coarse mesh), in our case $V_H$, and a fine part, in our case $\Vf$. The weak form is also decoupled into a coarse and a fine part. The method reads:
 find $\bar{u}\in V_H$ and $u'\in \Vf$ such that,
 \begin{gather*}
 a(\bar{u},\bar{v})+a(u',\bar{v})=G(\bar{v})\quad\text{for all } \bar{v}\in V_H,\\
 a(u',v')=G(v')-a(\bar{u},v')\quad\text{for all } v'\in \Vf.
 \end{gather*}
The fine scale solution is further decoupled over the coarse elements $T\in\triH$ and approximated using analytical techniques. Note that the fine scale solution $u'$ is an affine map of the coarse scale solution $\bar{u}$. If we let $u'\approx M\bar{u}+m$ and plug this in to the first equation we get a coarse stiffness matrix of the form $a(\bar{v}+M\bar{v},\bar{w})$, i.e., a non-symmetric bilinear form for a symmetric problem.

\subsection{The Multiscale Finite Element Method}
In \cite{MR1455261} the multiscale finite element method was first introduced. Here modified multiscale basis functions are computed numerically on sub-grids on each coarse element individually. The basis functions fulfill: find $\phi_{x,T}\in H^1_0(T)$
$$
a(\lambda_x-\phi_{x,T},v)=0\quad\text{for all } v\in H^1_0(T) \text{ and for all } T\in\triH.
$$
Here homogeneous Dirichlet boundary conditions are used on the boundary of each element $T$, i.e., the local problems are totally decoupled. To get a more accurate method one can improve the boundary conditions using information from the data $A$. A larger domain can also be considered (this procedure is referred to as over-sampling), see e.g.~\cite{MR1455261}. Note that since the coarse scale basis functions are modified (both trial and test space) the resulting method is symmetric.

\subsection{The Adaptive Variational Multiscale Method}
The modified basis function construction given by equation (\ref{e:Tnodal}-\ref{e:basiscoarse}) was first introduced in a variational multiscale framework in \cite{MR2161713,MR2319044}. In these papers the Scott-Zhang interpolation was used in the analysis and nodal interpolation in the discrete setting for the numerical examples. The modified basis functions where only used for the trail functions but not for the test functions. A fine scale correction based on the right hand side data was also included. In \cite{MRX} the modified basis functions were used for both trail and test functions. The exponential decay of the modified basis functions, with respect to number of coarse layers of elements in the vertex patches, have been demonstrated numerically in all these works, see \cite{MR2553176,MRX}.

The adaptive variational multiscale method has been extended to convection dominated problems and problems in mixed form \cite{MRX}. A posteriori error bounds have been derived and adaptive algorithms designed where the local mesh and patch size are chosen automatically in order to reduce the error.

\subsection{Application of the Presented Analysis}
The convergence proof in this paper gives a valid bound also as $h\rightarrow 0$ independent of the patch size and coarse mesh size. The proof does not rely on regularity of the solution and gives a very explicit expression for the rate of convergence. The present analysis confirms the numerical results in \cite{MR2553176,MRX} and gives the symmetric version of the method, where both trail and test space are modified, the solid theoretical foundation it has previously been missing. The analysis also justifies the use of a posteriori error bounds for adaptivity \cite{MR2319044,MRX} because we can now prove that the quantities measured on the patch boundary decays exponentially in the number of coarse layers.

For the variational multiscale method this result says that it is important to allow larger subgrid patches than just one coarse element. This will result in overlap but the local problems are totally decoupled and we have in previous works demonstrated how adaptivity can be used to only solve local problems where it is needed, see for instance \cite{MR2319044,MRX}. For the multiscale finite element method the analysis is not directly applicable since the fine scale space $\Vf$ is not used. It is the decay in this space which we have proven to be exponential (in number of coarse layers of elements in the subgrid). If this decay is not present, inhomogeneous boundary conditions are instead needed for the subgrid problems. To the best of our knowledge, such constructions have only been proved to be accurate in special settings, e.g.,~periodic coefficients.

\bibliographystyle{amsplain}
\bibliography{multiscale}

\end{document}